%% file: IRK.tex
\RequirePackage[table]{xcolor}
\documentclass{article}
\usepackage{geometry}
\usepackage{pdflscape}
\usepackage{graphicx}
\usepackage{amsmath,amssymb}
\usepackage{bm}
\usepackage{booktabs}
\usepackage{array}
\usepackage{subcaption}
\usepackage{siunitx}
\usepackage{hyperref}
\usepackage{algorithm}
\usepackage[noend]{algpseudocode}
\usepackage{listings,xcolor}
\usepackage{varwidth}
\usepackage{authblk}

\title{Stage-parallel fully implicit Runge-Kutta solvers for discontinuous 
Galerkin fluid simulations}
\author[1]{Will Pazner}
\author[2]{Per-Olof Persson}
\affil[1]{Division of Applied Mathematics, Brown University, Providence, RI, 02912}
\affil[2]{Department of Mathematics, University of California, Berkeley, Berkeley, CA, 94720-3840}

\date{}

\newcommand{\keywords}[1]{ {\small \textbf{\textit{Keywords:}} #1} }

\begin{document}
\maketitle
\begin{abstract}
In this paper, we develop new techniques for solving the large, coupled linear
systems that arise from fully implicit Runge-Kutta methods. This method makes
use of the iterative preconditioned GMRES algorithm for solving the linear
systems, which has seen success for fluid flow problems and discontinuous
Galerkin discretizations. By transforming the resulting linear system of
equations, one can obtain a method which is much less computationally expensive
than the untransformed formulation, and which compares competitively with other
time-integration schemes, such as diagonally implicit Runge-Kutta (DIRK)
methods. We develop and test several ILU-based preconditioners effective for
these large systems. We additionally employ a parallel-in-time strategy to
compute the Runge-Kutta stages simultaneously. Numerical experiments are
performed on the Navier-Stokes equations using Euler vortex and 2D and 3D NACA
airfoil test cases in serial and in parallel settings. The fully implicit Radau
IIA Runge-Kutta methods compare favorably with equal-order DIRK methods in terms
of accuracy, number of GMRES iterations, number of matrix-vector
multiplications, and wall-clock time, for a wide range of time steps.
\end{abstract}

\section{Introduction}
The discontinuous Galerkin method, introduced in 1973 by Reed and Hill for the
neutron transport equation \cite{Reed:1973}, has seen in recent years increased
interest for fluid dynamics applications \cite{Peraire:2011}. The discontinuous
Galerkin method is a high-order finite element method suitable for use on
unstructured meshes with polynomials of arbitrarily high degree. For many fluid
flow problems, explicit time integration methods have the downside of
restrictive time step conditions, partly because of the use of high-degree
polynomials but more fundamentally because of the need to employ highly graded
and/or anisotropic elements for many realistic flow problems
\cite{Persson:2006}. Therefore, for many applications it is desirable to use an
implicit time integration scheme.

Implicit time integration methods for DG have been much studied.
Multi-step backward differentiation formulas (BDF) and single-step diagonally
implicit Runge-Kutta (DIRK) methods have been applied to
discontinuous Galerkin discretizations for fluid flow problems
\cite{Persson:2011LES, Persson:2006}.
Nigro et.\ al have seen success applying
multi-stage, multi-step modified extended BDF (MEBDF) and two implicit advanced 
step-point (TIAS) schemes to the compressible Euler and Navier-Stokes equations
\cite{Nigro2014MEBDF, Nigro2014TIAS}.
Additionally, in \cite{Bassi2015}, Bassi
et al.\ have used linearly implicit Rosenbrock-type to integrate DG
discretizations for various fluid flow problems. The BDF and DIRK methods have
some limitations: BDF schemes can be $A$-stable only up to second-order (the
famous second Dahlquist barrier) \cite{Dahlquist:1963}, a severe limitation when
used in conjunction with a high-order spatial discretization. On the other hand,
there exist high-order $A$-stable (and even $L$-stable) DIRK schemes, but these
methods have a low stage-order, often resulting in order reduction when applied
to stiff problems \cite{Frank:1985kv}.

The Radau IIA methods, one class of the so-called fully implicit Runge-Kutta
(IRK) methods, are high-order, $L$-stable, and have relatively high stage order.
Consequently, these methods suffer less from order reduction than the
corresponding DIRK methods when applied to stiff problems. Furthermore, these
methods require only a small number of stages $s$, with the order of accuracy
given by $p = 2s - 1$. These methods have the drawback that each step involves
the solution of large, coupled linear systems of equations. The difficulty in
efficiently implementing such methods has caused them to remain not widely used
or studied for practical applications \cite{Carpenter:2005, Carpenter2003:ef}.
There has been previous work on improving the efficiency of solving these large,
coupled systems. In \cite{Jay1999}, Jay and Braconnier develop a parallelizable
preconditioner for IRK methods by means of Hairer and Wanner's
$W$-transformation. In \cite{DeSwart1998}, De Swart et al.\ have developed a
parallel software package for the four-stage Radau IIA method, and Burrage et
at.\ have developed a matrix-free, parallel implementation of the fifth-order
Radau IIA method in \cite{Burrage1999}.

In this paper, we develop a new strategy for efficiently solving the resulting
large linear systems by means of the iterative preconditioned GMRES algorithm. A
simple transformation of the linear system results in a significant reduction of
the cost per GMRES iteration. Furthermore, the block ILU(0) preconditioner, used
successfully with implicit time-integrators for the discontinuous Galerkin
method in \cite{Persson:2008by}, proves to be effective also for these large
systems.  A shifted, uncoupled, block ILU(0) factorization is also found to be
an effective preconditioner, with the advantage of allowing parallelism in time
by computing the stage solutions simultaneously.

The structure of this paper is as follows. In Section \ref{sec:eqns}, we
describe the governing equations and DG spatial discretization. In Section
\ref{sec:time}, we discuss the time integration schemes used in this paper.
Then, in Section \ref{sec:formulation}, we introduce the transformation used to
reduce the solution cost, and discuss the preconditioners used for the GMRES
method. Finally, in Section \ref{sec:numerical}, we perform numerical
experiments on a variety of test cases, in two and three spatial dimensions.

\section{Equations and spatial discretization}\label{sec:eqns}
The equations considered are the time-dependent, compressible Navier-Stokes
equations,
\begin{gather}
    \label{eq:ns-1}
    \frac{\partial\rho}{\partial t}+
        \frac{\partial}{\partial x_j}(\rho u_j)=0 \\
    \label{eq:ns-2}
    \frac{\partial}{\partial t}(\rho u_i) +
        \frac{\partial}{\partial x_j}
        (\rho u_i u_j) + \frac{\partial p}{\partial x_i}
        = \frac{\partial \tau_{ij}}{\partial x_j}
        \qquad \text{for $i=1,2,3,$}\\
    \label{eq:ns-3}
    \frac{\partial}{\partial t}(\rho E) + 
      \frac{\partial}{\partial x_j}
        \left(u_j(\rho E + p) \right) = -\frac{\partial q_j}{\partial x_j}
        + \frac{\partial}{\partial x_j} (u_i \tau_{ij}),
\end{gather}
where $\rho$ is the density, $u_i$ is the $i$th component of the velocity,
and $E$ is the total energy. The viscous stress tensor and heat flux are
given by
\begin{equation}
\tau_{ij} = \mu\left( \frac{\partial u_i}{\partial x_j} +
                      \frac{\partial u_j}{\partial x_i} -
                      \frac{2}{3}
                      \frac{\partial u_k}{\partial x_k}  \delta_{ij}
                      \right)
\qquad\text{and}\qquad
q_j = - \frac{\mu}{\rm Pr} \frac{\partial}{\partial x_j} \left(
E + \frac{p}{\rho} - \frac{1}{2} u_k u_k \right),
\end{equation}
where $\mu$ is the viscosity coefficient, and ${\rm Pr} = 0.72$ is the Prandtl
number. For an ideal gas, the pressure $p$ is given by the equation of state
\begin{equation} \label{eq:eos}
    p = (\gamma - 1)\rho \left( E - \frac{1}{2} u_k u_k \right),
\end{equation}
where $\gamma$ is the adiabatic gas constant. For the viscous problems, we
introduce an isentropic assumption of the form $p = K \rho^\gamma,$ for a given
constant $K$, as described in \cite{Kanner2015}. This additional simplification
can be thought of as an artificial compressibility model for the incompressible
flows simulated in Sections \ref{sec:naca-2d} and \ref{sec:naca-3d}. We prefer
to solve the compressible equations because they result in a system of ODEs
rather than differential-algebraic equations, and thus do not need specialized
projection-type solvers. This model decouples equation \eqref{eq:ns-3} from
equations \eqref{eq:ns-1} and \eqref{eq:ns-2}, and therefore results in one
fewer component to solve for. We remark that this simplification does not result
in a significant difference in the relative performance of the time integrators
studied in this paper, as shown in Section \eqref{sec:ev}, where the full
compressible Euler equations are solved, with no isentropic assumption.

We rewrite equations \eqref{eq:ns-1}, \eqref{eq:ns-2}, and \eqref{eq:ns-3} in
the form
\begin{equation}
    \label{eq:ns-conservation}
    \frac{\partial u}{\partial t} + \nabla \cdot \bm{F}_i(u)
        - \nabla \cdot \bm{F}_v (u, \nabla u) = 0,
\end{equation}
where $u$ is a vector of the conserved variables, and $\bm{F}_i, \bm{F}_v$ are
the inviscid and viscous flux functions, respectively. The spatial domain
$\Omega$ is discretized into a triangulation, and the solution $u$ is
approximated by piecewise polynomials of a given degree. Equation
\eqref{eq:ns-conservation} is discretized by means of the discontinuous Galerkin
method, where the viscous terms are treated using the compact DG (CDG) scheme
\cite{Peraire:2008}.

Using a nodal basis function for each component, we write the global solution 
vector as $\bm{u}$, and obtain a semi-discrete system of ordinary differential 
equations of the form
\begin{equation}
   \bm{M} \frac{\partial \bm{u}}{\partial t} = \bm{f}(\bm{u}),
\end{equation}
where $\bm{M}$ is the mass matrix, and $\bm{f}$ is a nonlinear function of the
$n$ unknowns $\bm{u}$. The standard method of lines approach allows for the
solution of this system of ordinary differential equations by means of a range
of numerical time integrators, such as the implicit Runge-Kutta methods that
are the focus of this paper.

\subsection{Block structure of the Jacobian}
We consider the vector of unknowns $\bm{u}$ to be ordered such that the $m$
degrees of freedom associated with one element of the triangulation appear
consecutively. We suppose that there are a total of $T$ elements, such that
there are a total of $n = Tm$ degrees of freedom. Then, the Jacobian matrix
$\bm{J}$ can be seen as a $T \times T$ block matrix, with blocks of size $m
\times m$. The $i$th row consists of blocks on the diagonal, and in columns $j$,
where elements $i$ and $j$ share a common edge, such that the total number of
off-diagonal blocks in the $i$th row is equal to the number of neighbors of
element $i$. We note that the off-diagonal blocks of size $m \times m$ are
themselves sparse, but for the sake of simplicity we will consider them as dense
matrices. The mass matrix $\bm{M}$ is a $T \times T$ block diagonal matrix, with
blocks of size $m \times m$, and therefore matrices of the form $\alpha \bm{M} -
\beta \bm{J}$ have the same sparsity pattern as the Jacobian.

\section{Time integration}\label{sec:time}
In this paper, we will focus on the one-step, multi-stage
Runge-Kutta methods. Given initial conditions $\bm{u}_0 = \bm{u}(t_0)$, a
general $s$-stage, $p$th-order Runge-Kutta method for advancing the solution to
$\bm{u}_1 = \bm{u}(t_0 + \Delta t) + \mathcal{O}(\Delta t^{p+1})$ can be written
as
\begin{align}
   \label{eq:DG-RK-1}
   \bm{M} \bm{k}_i &= \bm{f}\left( t_0 + \Delta t c_i,
                  \bm{u}_0 + \Delta t \sum_{j=1}^s a_{ij} \bm{k}_j \right), \\
   \label{eq:DG-RK-2}
   \bm{u}_1 &= \bm{u}_0 + \Delta t \sum_{i=1}^s b_i \bm{k}_i,
\end{align}
where the coefficients $a_{ij}, b_i,$ and $c_i$ can be expressed compactly in
the form of the \textit{Butcher tableau},
\[
   \begin{array}{c|ccc}
       c_1    & a_{11} & \cdots & a_{1s} \\
       \vdots & \vdots & \ddots & \vdots \\
       c_s    & a_{s1} & \cdots & a_{ss} \\
       \hline
              & b_1    & \cdots & b_s
   \end{array}
   =
   {\setlength\extrarowheight{5pt}
   \begin{array}{c|c}
       \bm{c} & \bm{A}\\
       \hline
       & \bm{b}^T
   \end{array}}.
\]

If the matrix of coefficients $\bm{A}$ is strictly lower-triangular, then each
stage $\bm{k}_i$ only depends on the preceding stages, and the method is called
an \textit{explicit} Runge-Kutta method. In this case, each stage may be
computed by simply evaluating the function $\bm{f}$. If $\bm{A}$ is not strictly
lower-triangular, the method is called an \textit{implicit} Runge-Kutta method
(IRK). A particular class of implicit Runge-Kutta methods is those for which the
matrix $\bm{A}$ is lower-triangular. Such methods are called
\textit{diagonally-implicit} Runge-Kutta (DIRK) methods \cite{Alexander:1977dk}.
Implicit Runge-Kutta methods enjoy high accuracy and very favorable stability
properties, but computing the stages requires the solution of (in general
nonlinear) systems of equations. In the case of DIRK methods, since $\bm{A}$ is
lower-triangular, each stage $\bm{k}_i$ depends only on those stages $\bm{k}_j$,
$j \leq i$, requiring the sequential solution of $s$ systems, each of size $n$.
In contrast, general IRK methods couple all of the stages, resulting in one
nonlinear system of equations of size $s\times n$.

For the solution of stiff systems of equations, we are interested in those
methods that are $L$-stable, meaning that their stability region includes the
entire left half-place ($A$-stability), together with the additional criterion
that the stability function $R(z)$ satisfies $\lim_{z\to\infty} R(z) = 0$. In
the present study, we compare the efficiency and effectiveness of several
$L$-stable IRK and DIRK schemes. The methods considered are listed in Table
\ref{tab:schemes}. The IRK schemes considered are the Radau IIA schemes, which
are $L$-stable, $s$-stage schemes of order $2s - 1$ based on the Radau right
quadrature. The construction of these schemes can be found in
\cite{Hairer:2013tj}. The two-stage and three-stage Radau IIA methods are listed
as RADAU23 and RADAU35, respectively. The DIRK schemes considered are the
three-stage, third-order $L$-stable scheme denoted DIRK33, which is derived in
detail in \cite{Alexander:1977dk}, and the six-stage, fifth-order scheme
constructed in \cite{Boom:2013}, and denoted ESDIRK65. The latter scheme is an
\textit{explicit singly diagonally implicit} Runge-Kutta (ESDIRK) method,
meaning that the first diagonal entry of the Butcher matrix is zero, and the
remaining diagonally entries are nonzero and equal. In addition to the third- 
and fifth-order methods, we also consider the seventh- and ninth-order Radau IIA
methods, although we do not compare these methods to equal-order DIRK methods.
The Butcher tableaux for the methods considered are given in Appendix 
\ref{app:tableaux}.

\begin{table}[h!]
\centering
\caption{$L$-stable implicit Runge-Kutta schemes considered}
\label{tab:schemes}
\begin{tabular}{lccccc}
\toprule
Scheme         & Order & Total stages & Implicit stages  & Stage order &
  Leading error coefficient\\
\midrule
RADAU23        & 3         & 2        & 2   & 2   & $1.39 \times 10^{-2}$\\
DIRK33         & 3         & 3        & 3   & 1   & $2.59 \times 10^{-2}$\\
RADAU35        & 5         & 3        & 3   & 3   & $1.39 \times 10^{-4}$\\
ESDIRK65       & 5         & 6        & 5   & 2   & $5.30 \times 10^{-4}$\\
RADAU47        & 7         & 4        & 4   & 4   & $7.09 \times 10^{-7}$\\
RADAU59        & 9         & 5        & 5   & 5   & $2.19 \times 10^{-9}$\\
\bottomrule
\end{tabular}
\end{table}

Also of interest is the phenomenon of \textit{order reduction}, whereby, when
applied to stiff problems,  the overall order of accuracy is reduced from $p$ to
the \textit{stage order} of the method (denoted $q$) \cite{Frank:1985kv}. The
stage order $q$ is defined as $q = \min\{ p, q_i \}$, for $i=1,\ldots,s$, where
$q_i$ is defined by
\begin{equation}
   \bm{u}(t_0 + \Delta t c_i) = \bm{u}_0 + \Delta t \sum_{j=1}^s a_{ij}
    \bm{k}_j + \mathcal{O}(\Delta t^{q_i + 1}).
\end{equation}
It can be shown that the maximum stage order for any DIRK method is 2
\cite{Hairer:2013vl}, whereas the stage order for the Radau IIA methods is given
by the number of stages, $q=s$ \cite{Lambert:1991vi}. The DIRK33 has stage order
of $q=1$. An advantage of the ESDIRK methods such as ESDIRK65 is that they have
stage order of $q=2$ \cite{Kvaerno:2004}.

The Radau IIA methods are very attractive because of their high order of
accuracy, small number of stages, high stage order, and $L$-stability, but
solving the coupled system of $s \times n$ equations is computationally
expensive. Supposing that we solve the nonlinear system of equations for the
stages $\bm{k}_i$ by means of Newton's method, then at each iteration we must
solve a linear system of equations by inverting the Jacobian matrix of the
right-hand side, $\bm{f}(t, \bm{u})$. Assuming a dense Jacobian matrix, and
solution via Gaussian elimination (or LU factorization), then the cost of
performing a linear solve scales as the cube of the number of unknowns.
Therefore, the cost per linear solve for a general IRK method is
$\mathcal{O}(s^3n^3)$, whereas the cost per solve for a DIRK method scales like
$\mathcal{O}(sn^3)$. In \cite{Butcher:1976tt}, Butcher describes how to
transform the resulting set of linear equations to reduce the computational work
for solving the IRK systems to $\mathcal{O}(2sn^3)$.  Despite this reduction in
computational complexity, the cost of solving the large systems of equations has
proven in practice to be prohibitive \cite{Carpenter:2005}. On the other hand,
DIRK methods have proven be popular and effective for solving computational
fluid dynamics problems \cite{Bijl:2002}, at the cost of lower stage order and
an increased number of stages.


\section{Efficient solution of implicit Runge-Kutta systems}
\label{sec:formulation}
In this section we describe a method for efficiently solving the systems arising
from general IRK methods when applied to discontinuous Galerkin discretizations.
The Jacobian matrices of the function $\bm{f}$ are sparse, block-structured
matrices, which lend themselves to solution via iterative Krylov subspace
methods. In particular, we consider the solution of these systems by means of
the GMRES method with a zero fill-in block ILU(0) preconditioner, as in
\cite{Persson:2008by}. In this case, each iteration of the GMRES method requires
one matrix-vector multiplication, and one application of the ILU(0)
preconditioner. In order to efficiently solve the linear systems resulting from
IRK methods, we will rewrite the system of equations in such a way so as to
reduce the cost of a matrix-vector multiplication from $s^2n^2$ to order
$s n^2$.

\subsection{Transformation of the system of equations}
Recalling equation that the stages $\bm{k}_i$ are given by the equation
\begin{equation}
   \label{eq:IRK-stage-eqn}
   \bm{M} \bm{k}_i = \bm{f}\left( t_0 + \Delta t c_i,
                  \bm{u}_0 + \Delta t \sum_{j=1}^s a_{ij} \bm{k}_j \right),
\end{equation}
we define $\bm{K}$ to be the concatenation of the vectors $\bm{k}_i$, $\bm{U}_0$
to be the concatenation of $s$ copies of $\bm{u}_0$, and $\bm{F}$ the function
$\bm{f}$ applied component-wise on these vectors. Then, we rewrite
equation \eqref{eq:IRK-stage-eqn} in vector form as
\begin{equation}
   (\bm{I}_s \otimes \bm{M})\bm{K} =
       \bm{F}\left(t_0 + \Delta t \bm{c},
                   \bm{U}_0 + \Delta t (\bm{A} \otimes \bm{I}_n) \bm{K} \right),
\end{equation}
where $\otimes$ is the Kronecker product, and $\bm{I}_n$ is the $n\times n$
identity matrix.

This nonlinear system can be solved by means of Newton's method, which
will require solving at each step a linear system of the form
\begin{equation}
\label{eq:newton-system}
\left(\left(\begin{array}{ccc}
   \bm{M} &        & 0 \\
          & \ddots &  \\
   0      &        & \bm{M}
\end{array}\right) - \Delta t
\left(\begin{array}{ccc}
   a_{11} \bm{J}_1 & \cdots & a_{1s} \bm{J}_1 \\
   \vdots & \ddots & \vdots \\
   a_{s1} \bm{J}_s & \cdots & a_{ss} \bm{J}_s,
\end{array}\right)\right)
\left(\begin{array}{ccc}
   \bm{k}_1 \\ \vdots \\ \bm{k}_s
\end{array} \right)
=
\left(\begin{array}{ccc}
   \bm{r}_1 \\ \vdots \\ \bm{r}_s
\end{array} \right)
\end{equation}
for the residual vectors $(\bm{r}_1, \ldots, \bm{r}_s)^T$, where the matrices on
the left-hand side are $s \times s$ block matrices blocks, with each block of
size $n \times n$. We use the following notation for the Jacobian matrix of
$\bm{f}$,
\[
   \bm{J}_i = \bm{J}_{\bm{f}} \left(t_0 + \Delta t c_i, \bm{u_0}
       + \Delta t\sum_{j=1}^s a_{ij} \bm{k}_j\right).
\]
We can rewrite equation \eqref{eq:newton-system} in the following form,
\begin{equation}
   \left( \bm{I}_s \otimes \bm{M} - \Delta t
       \left(
       \begin{array}{ccc}
       a_{11}J_1 & \cdots & a_{1s}J_1 \\
       \vdots & \ddots & \vdots \\
       a_{s1}J_s & \cdots & a_{ss}J_s
       \end{array}
       \right)
   \right) \bm{K} = \bm{R},
   \label{eq:newton-system-2}
\end{equation}
The sparsity pattern of the matrix in \eqref{eq:newton-system-2} is simply
that of the Jacobian matrix $J$ repeated $s \times s$ times, and we can conclude
that the cost of computing a matrix-vector product with this matrix is $s^2$ 
times that of computing the matrix-vector product of one Jacobian matrix.

In order to reduce the cost of the matrix-vector multiplication, we perform a
simple transformation in to rewrite \eqref{eq:newton-system-2} in a slightly
modified form. We begin by defining
\[
   \bm{w}_i = \sum_{j=1}^s a_{ij} \bm{k}_j,
\]
and similarly, we let $\bm{W}$ denote the vectors $\bm{w}_i$ stacked, such that
\[
   \bm{W} = (\bm{A} \otimes \bm{I}_n) \bm{K}.
\]
Then, we rewrite the nonlinear system of equations \eqref{eq:IRK-stage-eqn} in
terms of the variables $\bm{w}_i$ to obtain
\[
   \bm{M}\bm{k_i} = \bm{f}(t_0 + \Delta t c_i, \bm{u}_0 + \Delta t \bm{w}_i ),
\]
or, equivalently, in the case that the Butcher matrix $A$ is invertible,
\begin{equation} \label{eq:modified-IRK}
   (\bm{A}^{-1} \otimes \bm{M}) \bm{W} = \bm{F}(t_0 + \Delta t \bm{c},
       \bm{U}_0 + \Delta t \bm{W} ).
\end{equation}
In the transformed variables, the new solution $\bm{u}_1$ can be written as
\[
  \bm{u}_1 = \bm{u}_0 + \Delta t ( \bm{b}^T \bm{A}^{-1} \otimes \bm{I}_n)\bm{W},
\]
In the case of the Radau IIA methods, $\bm{b}^T \bm{A}^{-1} = (0,\ldots,0,1)$, 
and so this further simplifies to
\[
  \bm{u}_1 = \bm{u}_0 + \Delta t \bm{w}_s.
\]
Solving equation \eqref{eq:modified-IRK} with Newton's method gives rise to the 
linear system of equations
\begin{equation}
   \left( A^{-1} \otimes \bm{M} - \Delta t
   \left( \begin{array}{cccc}
       \bm{J}_1    & 0           & \cdots & 0 \\
       0           & \bm{J}_2    & \cdots & 0 \\
       \vdots      & \vdots      & \ddots & \vdots \\
       0           & 0           & \cdots & \bm{J}_s
   \end{array} \right)
   \right) \bm{W} = \bm{R}.
\end{equation}
The advantage of this formulation is that the resulting system enjoys greater 
sparsity. The resulting matrix is a $s \times s$ block matrix, with multiples of
the mass matrix in every off-diagonal block, and with matrices of the form
$(\bm{A}^{-1})_{ii} \bm{M} - \Delta t \bm{J}_i$ along the diagonal. Computing a
matrix-vector product with this $s \times s$ block matrix
requires performing $s$ matrix-vector multiplications with the mass matrix, and 
$s$ matrix-vector multiplications with a Jacobian $\bm{J}_i$. 
Therefore the cost of computing such products scales as $s$ times the cost of 
computing one matrix-vector product with the Jacobian matrix.

We additionally remark that the fully-implicit IRK methods 
requiring storing each of the $s$ Jacobian matrices $\bm{J}_i$, resulting in
memory usage that is $s$-times that of the DIRK methods. A further advantage of
the transformed system of equations is that the memory requirements for the
ILU-based preconditioners are reduced, as discussed in the following sections.

\subsection{Preconditioning}
The use of an appropriate preconditioner is essential in
accelerating the convergence of a Krylov subspace method such as GMRES. We
briefly describe the block ILU(0) preconditioner from \cite{Persson:2008by}.

\subsubsection{Block ILU(0) preconditioner}
The block ILU(0) (or zero fill-in) preconditioner is a method for obtaining
block lower- and upper-triangular matrices $\tilde{\bm{L}}$ and $\tilde{\bm{U}}$
given a block sparse matrix $\bm{B}$. These matrices are obtaining by computing
the standard block LU factorization, but discarding any blocks which do not
appear in the sparsity pattern of $\bm{B}$. We denote the block in the $i$th row
and $j$th column as $\bm{B}_{ij}$. The ILU(0) algorithm can be written as shown
in Algorithm \ref{alg:ilu}.
\begin{algorithm}
\caption{Block ILU(0) algorithm}
\label{alg:ilu}
\begin{algorithmic}[1]
   \For{$i=1$ to $T$}
       \For{neighbors $j$ of $i$ with $j > i$}
           \State $\bm{B}_{ji} \gets \bm{B}_{ji} \bm{B}_{ii}^{-1}$
           \State $\bm{B}_{jj} \gets \bm{B}_{jj} - \bm{B}_{ji}\bm{B}_{ij}$
           \For{neighbors $k$ of $j$ and $i$ with $k > j$}
               \State $\bm{B}_{jk} \gets \bm{B}_{jk} - \bm{B}_{ji}\bm{B}_{ik}$
           \EndFor
       \EndFor
   \EndFor
   \State $\tilde{\bm{L}} \gets I+\text{strict block lower triangle of } \bm{B}$
   \State $\tilde{\bm{U}} \gets \text{block upper triangle of } \bm{B}$
\end{algorithmic}
\end{algorithm}

If we impose the condition on the triangulation of our domain that, if elements
$j$ and $k$ both neighbor element $i$, then elements $j$ and $k$
are not neighbors of each other, then the ILU(0) algorithm has the particularly
simple form, show in Algorithm \ref{alg:simp-ilu}. In practice, most well-shaped
meshes satisfy this condition and henceforth we will use this simpler algorithm.
\begin{algorithm}
\caption{Simplified block ILU(0) algorithm}
\label{alg:simp-ilu}
\begin{algorithmic}[1]
   \For{$i=1$ to $T$}
       \For{neighbors $j$ of $i$ with $j > i$}
           \State $\bm{B}_{ji} \gets \bm{B}_{ji} \bm{B}_{ii}^{-1}$
           \State $\bm{B}_{jj} \gets \bm{B}_{jj} - \bm{B}_{ji}\bm{B}_{ij}$
       \EndFor
   \EndFor
   \State $\tilde{\bm{L}} \gets I+\text{strict block lower triangle of } \bm{B}$
   \State $\tilde{\bm{U}} \gets \text{block upper triangle of } \bm{B}$
\end{algorithmic}
\end{algorithm}

\subsection{Preconditioning the large system} \label{sec:precond}
In the case of the general IRK methods, we must solve systems of the form
\begin{equation}
   \label{eq:big-irk}
   \bm{B}\bm{W} = \bm{R}, \qquad \bm{B} = 
   \left( \bm{A}^{-1} \otimes \bm{M} - \Delta t
   \left( \begin{array}{cccc}
       \bm{J}_1    & 0           & \cdots & 0 \\
       0           & \bm{J}_2    & \cdots & 0 \\
       \vdots      & \vdots      & \ddots & \vdots \\
       0           & 0           & \cdots & \bm{J}_s
   \end{array} \right)\right).
\end{equation}
We remark that the matrix $\bm{B}$ can now be considered as a $s \times s$ block
matrix, with blocks of size $n \times n$. Each $n \times n$ block is itself  a
$T \times T$ block matrix, with subblocks of size $m \times m$. We introduce the
notation $\bm{B}_{k\ell,ij}$ to denote the $(i, j)$ subblock of the $(k,\ell)$
block of $\bm{B}$. That is to say, $\bm{B}_{k\ell,ij}$ is the $(i, j)$ block of
the matrix $(\bm{A})^{-1}_{k\ell}\bm{M} - \delta_{k\ell}\Delta t \bm{J}_k$,
where $\delta_{k\ell}$ is the Kronecker delta.

\subsubsection{Stage-coupled block ILU(0) preconditioner}
We consider two preconditioners for this large $sn \times sn$ system. The first
is the standard (stage-coupled) block ILU(0) preconditioner, which can be
computed using Algorithm \ref{alg:big-ilu}. We note that this preconditioner
couples all $s$ stages of the method. This preconditioner requires
storing $s$ Jacobian-sized diagonal blocks, and $s^2 - s$ off-diagonal blocks, 
each the same size as the mass matrix.
\begin{algorithm}
\caption{ILU(0) algorithm for IRK systems of the form \eqref{eq:big-irk}}
\label{alg:big-ilu}
\begin{algorithmic}[1]
   \For{$k = 1$ to $s$}
       \For{$i=1$ to $T$}
           \For{neighbors $j$ of $i$ with $j > i$}
               \State $\bm{B}_{kk,ji} \gets \bm{B}_{kk,ji} \bm{B}_{kk,ii}^{-1}$
               \State $\bm{B}_{kk,jj} \gets \bm{B}_{kk,jj}
                   - \bm{B}_{kk,ji}\bm{B}_{kk,ij}$
           \EndFor
           \For{$\ell = k+1$ to $s$}
               \State $\bm{B}_{\ell k,ii} \gets \bm{B}_{\ell k,ii}
                   \bm{B}_{kk,ii}^{-1}$
               \State $\bm{B}_{\ell \ell,ii} \gets \bm{B}_{\ell \ell,ii}
                   - \bm{B}_{\ell k,ii} \bm{B}_{kk,ii}$
           \EndFor
       \EndFor
   \EndFor
   \State $\tilde{\bm{L}} \gets I+\text{strict block lower triangle of } \bm{B}$
   \State $\tilde{\bm{U}} \gets \text{block upper triangle of } \bm{B}$
\end{algorithmic}
\end{algorithm}

\subsubsection{Stage-uncoupled, shifted ILU(0) preconditioner}
In order to avoid the above coupling of the stages, we can compute a simplified
preconditioner in the form of the following block matrix,
\begin{equation}
   \left( \begin{array}{cccc}
   \tilde{\bm{L}}_1 \tilde{\bm{U}}_1 & 0 & \cdots & 0 \\
   0 & \tilde{\bm{L}}_2 \tilde{\bm{U}}_2 & \cdots & 0 \\
   \vdots & \vdots & \ddots & \vdots \\
   0 & 0 & \cdots & \tilde{\bm{L}}_s \tilde{\bm{U}}_s
   \end{array} \right),
\end{equation}
where $\tilde{\bm{L}}_i \tilde{\bm{U}}_i$ is the block ILU(0) factorization of a
matrix of the form $\left(A^{-1}_{ii} + \alpha_i\right) \bm{M} - \Delta t
\bm{J}_i$. We let $\alpha_i$ denote a \textit{shift}, so that the standard
unshifted factorization corresponds to $\alpha_i = 0$. The so-called shifted ILU
factorization, described by Manteuffel in \cite{Manteuffel1980}, can result in
eigenvalues clustered away from the origin, and hence faster convergence in
GMRES. Indeed, our experience shows that the unshifted preconditioner
underperforms certain other choices of shift.

The preconditioner has several advantages over the stage-coupled block ILU(0)
preconditioner. The first is that it is easily constructed using an
already-implemented block ILU(0) factorization of the Jacobian matrix. The
second is that none of the stages are coupled, allowing for both efficient
computation and application. In particular, this has implications for the
parallelization of the preconditioner, which we discuss in Section
\ref{sec:parallel}. Finally, as this preconditioner does not 
include any off diagonal blocks, the memory requirements are exactly $s$ times 
that of the standard block ILU(0) used for the DIRK methods.

As mentioned, the unshifted preconditioner, with $\alpha_i = 0$ for all $i$,
such that $\tilde{\bm{L}}_i\tilde{\bm{U}}_i$ is the ILU(0) factorization of the
$i$th diagonal block of the matrix $\bm{B}$, is a natural choice. This choice of
coefficients ignores all the off-diagonal mass matrices. By making certain
judicious choices of the coefficients $\alpha_i$, we can attempt to compensate
for the off-diagonal blocks by adding multiples of the mass matrix back
to the diagonal entries. In particular, our numerical experiments have shown
that setting $\alpha_i = \sum_{j \neq i} \left| A^{-1}_{ji} \right|$ results in
a more efficient preconditioner, requiring fewer GMRES iterations in order to
converge to a given desired tolerance.

\subsection{Computational cost and memory requirements}
In order to compare the computational cost of the transformed IRK implementation
described above with both that of the untransformed formulation, and with the
usual DIRK methods, we summarize the computational cost associated with solving
the resulting linear systems. We note that the IRK methods require the solution
of one large, coupled system, whereas the DIRK methods require the solution of
$s$ smaller systems. In Table \ref{tab:irk-iter-costs} we record the leading
terms of the computational cost of operations required to be performed every
iteration. We recall that $s$ is the number of Runge-Kutta stages, $m$ is the
number of degrees of freedom per mesh element, $T$ is the total number of
elements in the mesh, and $r$ is the number of neighbors per element. Here we
assume that the $m \times m$ blocks of the Jacobian matrix are dense, and hence
require $2m^2$ floating point operations per matrix-vector multiply. Computing
the preconditioner requires the LU factorization of the diagonal blocks,
which incurs a cost of $\frac{2}{3}m^3$ floating point operations per block.
Each iteration in Newton's method requires re-evaluation of the
Jacobian matrix, and therefore also the re-computation of the preconditioner.
Hence, the preconditioner must be computed once per linear solve. The costs
associated with these calculations are listed in Table
\ref{tab:irk-precond-cost}. 

\begin{table}[H]
\centering
\caption{Per-iteration computational costs for solving implicit Runge-Kutta
systems}
\label{tab:irk-iter-costs}
\begin{tabular}{ll}
\toprule
Operation & Cost (leading term) \\
\midrule
Untransformed IRK matrix-vector product  & $s^2 m^2 (r + 1) T $ \\
Transformed IRK matrix-vector product    & $s m^2 (r + s) T $ \\
DIRK matrix-vector product               & $m^2 (r + 1) T$\\
\midrule
Coupled preconditioner application (IRK)        & $s m^2 (r + s) T$\\
Uncoupled preconditioner application (IRK)      & $s m^2 (r + 1) T$\\
Preconditioner application (DIRK)               & $m^2 (r + 1) T$\\
\bottomrule
\end{tabular}
\end{table}

\begin{table}[H]
\centering
\caption{Computational cost of computing the preconditioner (once per solve)
systems}
\label{tab:irk-precond-cost}
\begin{tabular}{ll}
\toprule
Operation & Cost (leading term) \\
\midrule
Computing coupled block ILU(0) preconditioner (IRK) & $s (m^3 + (r + s)m^2) T$\\
Computing uncoupled ILU(0) preconditioner (IRK)     & $s (m^3 + rm^2) T$\\
Computing block ILU(0) preconditioner (DIRK)        & $(m^3 + rm^2) T$\\
\bottomrule
\end{tabular}
\end{table}

We remark that each GMRES iteration using the formulation described in Section
\ref{sec:formulation} requires a factor of $s$ fewer floating-point operations
per iteration than the na\"ive IRK implementation. The stage-uncoupled IRK
preconditioner is less expensive to both compute and apply than the
stage-coupled preconditioner. We also note that for equal order of accuracy, the
Radau IIA IRK methods require fewer implicit stages than do the DIRK methods.
Each such implicit stage incurs the cost of assembling the Jacobian matrix. This
cost is problem-dependent, but is in general non-trivial, and in the model
problems considered in this paper, it scales like $\mathcal{O}(m^3 T)$.

Finally, we present the memory requirements for the IRK and DIRK
methods, and the stage-coupled and uncoupled preconditioners in Table
\ref{tab:irk-memory}. We note that for the transformed IRK methods, only the $s$
Jacobian matrices $\bm{J}_i$ need to be stored. The stage-coupled block ILU(0)
preconditioner requires an additional $s^2 - s$ off-diagonal blocks, which have
the same block-diagonal structure as the mass matrix, each having $m^2 T$
nonzero entries. The stage-uncoupled preconditioner does not require these 
off-diagonal blocks, and therefore its memory requirements are exactly $s$ 
times that of the DIRK block ILU(0) preconditioner. The block ILU(0) 
preconditioner for the untransformed system of equations would require storing 
an additional $(s^2 - s) r m^2 T$ nonzero entries in the off-diagonal blocks.

\begin{table}[H]
\centering
\caption{Memory requirements for the Jacobian matrices and preconditioners}
\label{tab:irk-memory}
\begin{tabular}{ll}
\toprule
Method & Memory required \\
\midrule
IRK Jacobian matrix                       & $s(r+1)m^2 T$\\
DIRK Jacobian matrix                      & $(r+1)m^2 T$\\
Coupled block ILU(0) preconditioner (IRK) & $s m^2 (s+r) T$\\
Uncoupled ILU(0) preconditioner (IRK)     & $s m^2 (r+1) T$\\
Block ILU(0) preconditioner (DIRK)        & $m^2 (r+1) T$\\
\bottomrule
\end{tabular}
\end{table}

\subsection{Stage-parallelism and partitioned ILU}\label{sec:parallel}
In order to parallelize the computations, the spatial domain is decomposed into
several subdomains. The compact stencil of the CDG scheme allows for very low
communication costs between processes for residual evaluation and Jacobian
assembly operations. The matrix-vector multiplications, which constitute the
bulk of the work for the linear solve, also scale well in terms of communication
for the same reason. In order to parallelize the ILU(0) preconditioner, the
contributions between elements in different partitions are ignored, allowing
each process to compute the ILU(0) factorization independently. Because these
contributions are ignored, we find that the number of GMRES iterations required
to converge increases as the number of domain partitions increases
\cite{Persson2009scalable}. In fact, it is easy to see that when the number of
partitions is equal to the number of mesh elements, the preconditioner simply
reduces to the block Jacobi preconditioner. However, in general the partitioned
ILU(0) preconditioner is found to be superior to the block Jacobi preconditioner.

These considerations apply equally to both the DIRK and general IRK methods,
using the stage-coupled ILU(0) preconditioner. If we use the stage-uncoupled
ILU(0) preconditioner, then we are able to decompose the domain into a factor of
$s$ fewer partitions. We then assume that the number of processes is equal to
the number of mesh partitions times the number of stages. The processes are
first divided into groups according to the mesh partitioning such that each
group consists of $s$ processes. Within each group, each process is then
assigned to one stage of the IRK method. Thus, when assembling the block matrix
of the form \eqref{eq:big-irk}, the Jacobian matrices for all of the stages are
computed in parallel. This does not require any communication between the
groups. Then, since the stage-uncoupled block ILU(0) preconditioner does not
take into account any of the off-diagonal blocks, the preconditioner can also be
computed without any inter-stage communication. Similarly, each application of
the preconditioner can be computed in parallel over all the stages without any
communication. When computing matrix-vector products, the products with the
Jacobian blocks for each stage are computed in parallel, and the products with
the mass matrix blocks must be communicated within the stages. It is possible to
overlap the communication with the computation of the matrix-vector product with
the stage-Jacobian block, such that the cost of communication is negligible. The
main advantage of this parallelization scheme is that the mesh is decomposed
into a factor of $s$ fewer partitions. Therefore, the effect of ignoring the
coupling between regions in the ILU(0) factorization is lessened, and the result
is a more efficient preconditioner. In Sections \ref{sec:parallel-results} and
\ref{sec:naca-3d}, we numerically study the effects of parallelizing these
preconditioners.

\begin{figure}[H]
   \centering
   \includegraphics[width=5in]{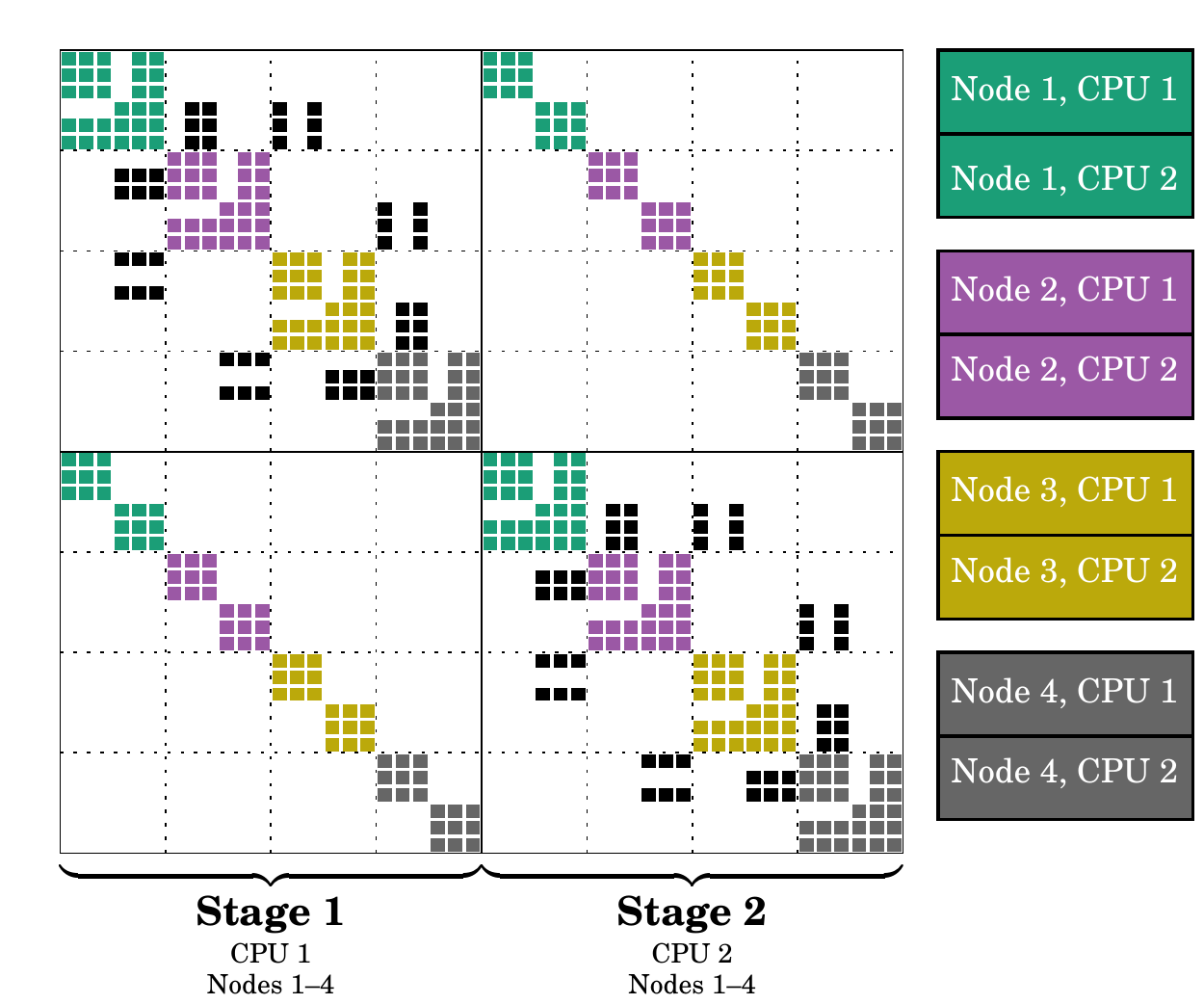}
   \caption[Stage-parallel partitioning of the Jacobian matrix]
           {Stage-parallel partitioning of the matrix
            $\bm{B} = \bm{A}^{-1} \otimes \bm{M} - \Delta t
            \left(
            \begin{smallmatrix}
            \bm{J}_1 & 0        \\
            0        & \bm{J}_2 \\
            \end{smallmatrix}
            \right)$.
            The color of the matrix entries indicates the node to which it
            belongs. The entries in the first block-column all belong to the
            first CPU of each node, and the entries the second block-column all
            belong to the second CPU of each node.}
   \label{fig:parallel-matrix}
\end{figure}

In the context of hybrid shared-distributed memory systems, it is possible to
further reduce the communication cost of the stage-parallel IRK algorithm. On
such a system, groups of compute cores called nodes have access to the same
shared memory. As a consequence, intranode communication is much faster than
internode communication. If each stage-group of $s$ processes described above
is located on one node, and none of the groups are split across nodes, then the
solutions are only communicated within a node, resulting in negligible internode
communication costs. An illustration of such an arrangement is shown in Figure
\ref{fig:parallel-matrix}. The example shown is a scalar problem on a mesh with
eight elements, decomposed into four partitions. The hypothetical architecture
consists of four compute nodes, each with two CPUs with shared memory. For a
two-stage IRK method, each partition of the mesh belongs to a different node,
and each of the stages for a given partition belong to different CPUs within one
node.

\section{Numerical results}\label{sec:numerical}
\subsection{Euler vortex} \label{sec:ev}

We solve the compressible Euler equations of gas dynamics, which are given by 
equations \eqref{eq:ns-1} through \eqref{eq:ns-3}, with the second-order terms 
removed. The equation of state is given by \eqref{eq:eos}.
%
\begin{figure*}[b!]
   \centering
   \begin{subfigure}[t]{0.5\textwidth}
       \centering
       \includegraphics[width=2.3in]{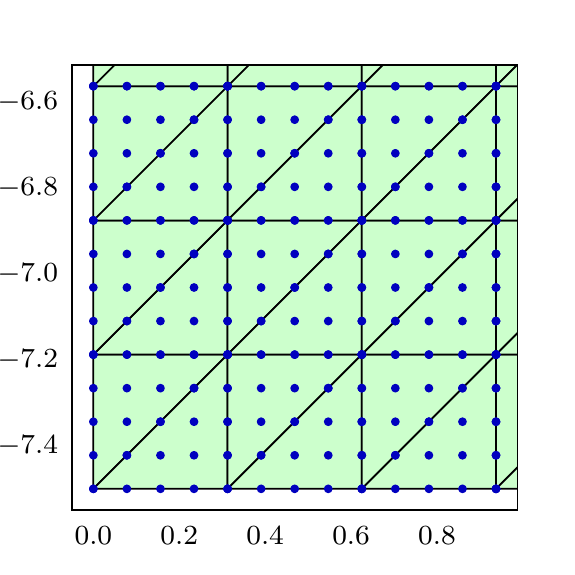}
       \caption{Bottom-left corner of mesh, with $p=4$ DG nodes}
       \label{fig:ev-mesh}
   \end{subfigure}%
   ~
   \begin{subfigure}[t]{0.5\textwidth}
       \centering
       \includegraphics[width=3in]{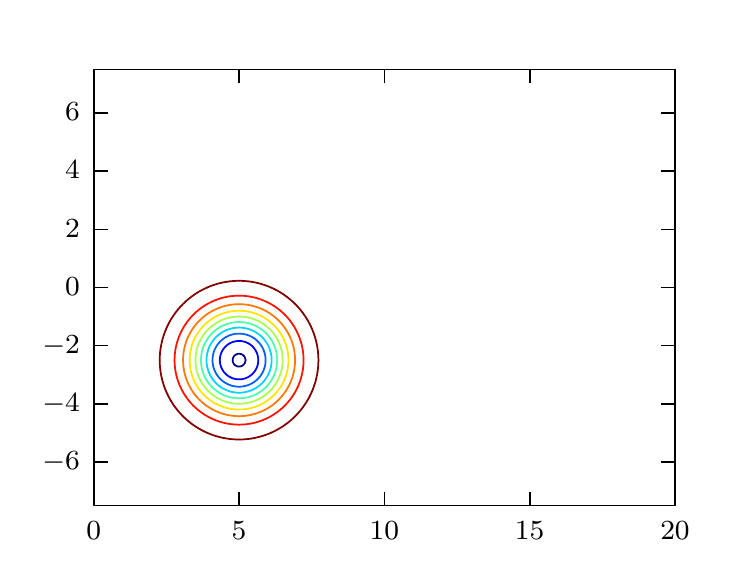}
       \caption{Density contours for the initial conditions}
       \label{fig:ev-ic}
   \end{subfigure}
   \caption{Compressible Euler vortex}
\end{figure*}
We consider the model problem of an unsteady compressible vortex in a
rectangular domain \cite{Wang:2013gj}. The domain is taken to be the rectangle $
[0, 20] \times [-7.5, 7.5]$, and the vortex is initially centered at $(x_0, y_0)
= (5, -2.5)$. The vortex is moving with the free-stream at an angle of $\theta$.
This problem is particularly useful as a benchmark because the exact solution
is given by the following analytic formulas, allowing for convenient computation
of the numerical accuracy. The exact solution at $(x, y, t)$ is given by
\begin{gather}
\label{eq:ev-exact-1}
   u = u_\infty \left( \cos(\theta) - \frac{\epsilon ((y-y_0)
       - \overline{v} t)}{2\pi r_c}
       \exp\left( \frac{f(x,y,t)}{2} \right) \right),\\
\label{eq:ev-exact-2}
   u = u_\infty \left( \sin(\theta) - \frac{\epsilon ((x-x_0)
       - \overline{u} t)}{2\pi r_c}
       \exp\left( \frac{f(x,y,t)}{2} \right) \right),\\
\label{eq:ev-exact-3}
   \rho = \rho_\infty \left( 1 -
       \frac{\epsilon^2 (\gamma - 1)M^2_\infty}{8\pi^2} \exp((f(x,y,t))
       \right)^{\frac{1}{\gamma-1}}, \\
\label{eq:ev-exact-4}
   p = p_\infty \left( 1 -
       \frac{\epsilon^2 (\gamma - 1)M^2_\infty}{8\pi^2} \exp((f(x,y,t))
       \right)^{\frac{\gamma}{\gamma-1}},
\end{gather}
where $f(x,y,t) = (1 - ((x-x_0) - \overline{u}t)^2 - ((y-y_0) -
\overline{v}t)^2)/r_c^2$, $M_\infty$ is the Mach number, $u_\infty,
\rho_\infty,$ and $p_\infty$ are the free-stream velocity, density, and
pressure, respectively. The free-stream velocity is given by $(\overline{u},
\overline{v}) = u_\infty (\cos(\theta), \sin(\theta))$. The strength of the
vortex is given by $\epsilon$, and its size is $r_c$.

For our test case, we choose parameters $\epsilon = 15$, $r_c = 1.5$, $M_0 =
0.5$, $\theta = \arctan(1/2)$. In order for the temporal error to dominate the
spatial error, we compute the solution on a fine mesh consisting of 6144 regular
right-triangular elements depicted in Figure \ref{fig:ev-mesh}. The finite
element space is chosen to be piecewise polynomials of degree 4, corresponding
to 15 nodes per element. The solution consists of four components, for a total
of 368,640 degrees of freedom. The initial conditions are shown in Figure
\ref{fig:ev-ic}.

\subsubsection{Comparison of preconditioners} \label{sec:ev-precond}
We first use this test case to compare the effectiveness of the preconditioners
discussed in Section \ref{sec:precond}. As a baseline, we will consider the DIRK
solver using the block ILU(0) preconditioner. We will then compare the
stage-coupled block ILU(0) and both shifted and unshifted stage-uncoupled ILU(0)
preconditioners for the IRK solver. In all cases, we require
Newton's method to converge to an infinity-norm tolerance of $10^{-8}$. We
compare the number of GMRES iterations required to converge to a
relative, preconditioned tolerance of $10^{-5}$. In order to make a
fair comparison between the methods, we will compute the number of $n \times n$
matrix-vector multiplications required per iteration of Newton's method, across
all the stages. We will refer to this quantity as the number of
\textit{equivalent multiplications}. For the DIRK methods, this number is equal
to the number of GMRES iterations times the number of implicit stages. For the
IRK methods, we recall that each multiplication by the large block matrix of the
form \eqref{eq:big-irk} essentially consists of $s$ $n \times n$
matrix-vector multiplications.

\begin{figure*}[t!]
   \centering
   \includegraphics[width=5.6in]{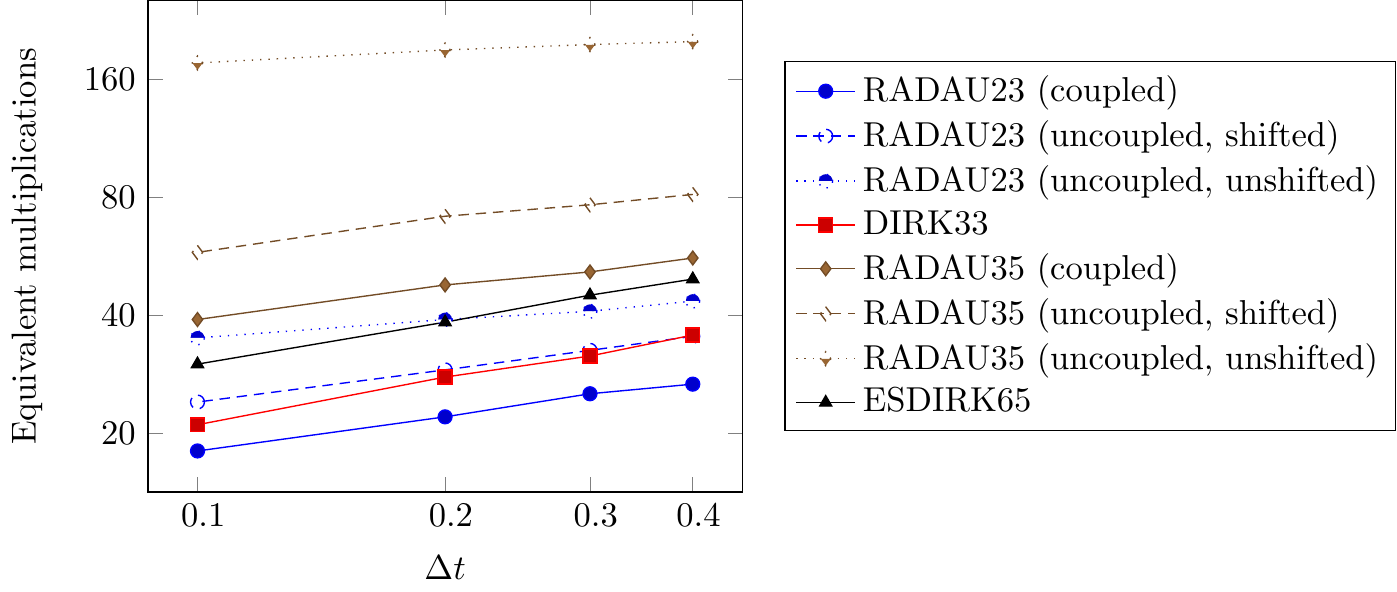}
   \caption{Log-log plot of number of average number of equivalent
            multiplications vs.\ $\Delta t$. Coupled preconditioners
            are shown in solid lines, uncoupled shifted preconditioners
            in dashed lines, and uncoupled, unshifted preconditioners
            in dotted lines.}
   \label{fig:ev-precond}
\end{figure*}

We compute 5 time steps in serial using representative time steps of $\Delta t =
0.4, 0.3, 0.2, $ and $0.1$. We then average the number of GMRES iterations
required per linear solve, and multiply this number by the number of implicit
stages. For the shifted stage-uncoupled block ILU(0) preconditioner, we choose a
shift of $\alpha_i = \sum_{j \neq i} \left| A^{-1}_{ji} \right|$, which has in
our experience resulted in the fastest convergence. We present the results in
the log-log plot shown in Figure \ref{fig:ev-precond}. We notice that the third
order Radau IIA method with the stage-coupled ILU(0) preconditioner requires
fewer matrix-vector multiplications than the corresponding third-order DIRK
method, while the stage-uncoupled, shifted ILU(0) preconditioner requires
roughly the same number of multiplications. For the fifth-order methods, both
the coupled and uncoupled preconditioners require more matrix-vector
multiplications than the corresponding ESDIRK method. For both third- and
fifth-order methods, the stage-uncoupled, unshifted preconditioner requires
greatly more matrix-vector multiplications than the other methods, and therefore
in the further test cases we will only consider the shifted preconditioner.

\subsubsection{Temporal accuracy}
Since the analytical solution to this test case is known, it is particularly
convenient to compare the accuracy of the time discretization schemes. The
solution is integrated until a final time of $t = \SI{60}{\s}$. For the
third-order methods, we choose time steps of $\Delta t = 0.4, 0.3, 0.2, 0.1,
0.075, 0.05, \SI{0.025}{\s}$. Because of the increased accuracy of the
fifth-order methods, we choose larger time steps of $\Delta t = 0.8, 0.6, 0.5,
0.4, 0.3, \SI{0.2}{\s}$ for the RADAU35 and ESDIRK65 methods. Time steps between
$0.5$ and $\SI{1.2}{\s}$ are chosen for the seventh- and ninth-order Radau IIA
methods. We approximate the $L^\infty$ error by comparing the numerical solution
with the known analytic solution at the DG nodes. We also estimate the order of
accuracy by comparing successive choices of $\Delta t$ and computing the rate of
convergence
\[
   r_i = \frac{ \log(L^\infty_{i+1} / L^\infty_i) }
               { \log(\Delta t_{i+1} / \Delta t_i) },
\]
where $L^\infty_i$ denotes the $L^\infty$ error of the numerical solution
computed using time step $\Delta t_i$. For each method, we present the
wall-clock time required to compute the solution in parallel on 16 cores. The
results for both the third-order and fifth-order methods are presented in Table
\ref{tab:ev-errors-3}. The theoretical order of accuracy is observed for all of
the methods used. Also listed is the ratio of the DIRK error to the Radau IIA
error.  Comparing the coefficients presented in Table
\ref{tab:schemes}, it can easily be shown that the leading coefficient of the
truncation error for the DIRK33 method is about 1.86 times larger than the
leading coefficient for the RADAU23 method. We see that the ratio of the errors
approaches this value as $\Delta t$ tends to zero. The leading coefficient of
the truncation error for the ESDIRK65 methods is about 3.82 times larger than
the leading coefficient for the RADAU35 method. In the fifth-order test case,
the ratio of the numerical errors is found to be closer to about 1.5, likely
because of the additional contribution of spatial discretization error.

\begin{table}[t!]
\footnotesize
\centering
\caption{$L^\infty$ error and runtime for the Euler vortex, DIRK and 
Radau IIA methods (wall-clock time presented for stage-coupled/uncoupled
preconditioners).}
\label{tab:ev-errors-3}
\begin{tabular}{l|lcc>{\hspace{1pc}}ccc>{\hspace{1pc}}c}
\toprule
          & \multicolumn{3}{c}{RADAU23} & \multicolumn{3}{c}{DIRK33} \\
$\Delta t$ & $L^\infty$ error & Order & Time C/UC (s)
          & $L^\infty$ error & Order & Time (s) & Ratio  \\
\midrule
0.2  &$7.2932\times 10^{-2}$& -  &820/829 &$1.0591\times 10^{-1}$&-   &1096& 1.452\\
0.1  &$9.7053\times 10^{-3}$&2.91&965/956 &$1.6115\times 10^{-2}$&2.72&1268& 1.660\\
0.075&$4.1135\times 10^{-3}$&2.98&1227/1223&$7.1881\times 10^{-3}$&2.81&1618& 1.747\\
0.05 &$1.2231\times 10^{-3}$&2.99&1734/1746&$2.2156\times 10^{-3}$&2.90&2291& 1.811\\
0.025&$1.7832\times 10^{-4}$&2.78&3214/3322&$3.0148\times 10^{-4}$&2.88&4212& 1.691\\
\toprule
          & \multicolumn{3}{c}{RADAU35} & \multicolumn{3}{c}{ESDIRK65} \\
$\Delta t$ & $L^\infty$ error & Order & Time C/UC (s)
          & $L^\infty$ error & Order & Time (s) & Ratio  \\
\midrule
0.6  &$5.8633\times 10^{-2}$& -  &843/775	&$8.3294\times 10^{-2}$&-   &697 & 1.421\\
0.5  &$2.4239\times 10^{-2}$&4.84&974/810	&$3.6075\times 10^{-2}$&4.59&802 & 1.488\\
0.4  &$7.8902\times 10^{-3}$&5.03&862/905	&$1.2605\times 10^{-2}$&4.71&947 & 1.598\\
0.3  &$1.8156\times 10^{-3}$&5.11&1079/1030&$3.0723\times 10^{-3}$&4.91&1197& 1.692\\
0.2  &$2.7294\times 10^{-4}$&4.67&1523/1478&$4.1279\times 10^{-4}$&4.95&1121& 1.512\\
\bottomrule
\end{tabular}
\begin{tabular}{l|lcc>{\hspace{1.5pc}}c|lcc}
\toprule
         & \multicolumn{3}{c}{RADAU47} & \multicolumn{4}{c}{RADAU59} \\
$\Delta t$ & $L^\infty$ error & Order & Time C/UC (s) &
$\Delta t$ & $L^\infty$ error & Order & Time C/UC (s)  \\
\midrule
1.0  &$1.8699\times 10^{-2}$& -  & 1222/1029  & 1.2 &$5.0425\times 10^{-3}$& -  & 2117/1533 \\
0.8  &$4.5413\times 10^{-3}$&6.34& 1158/1016  & 1.0 &$1.7294\times 10^{-3}$&5.87& 2094/1567 \\
0.6  &$8.7259\times 10^{-4}$&5.73& 1463/1307  & 0.8 &$4.6966\times 10^{-4}$&5.84& 1988/1554 \\
0.5  &$3.0265\times 10^{-4}$&5.81& 1732/1541  & 0.6 &$8.8679\times 10^{-5}$&5.79& 2544/2040 \\
\bottomrule
\end{tabular}
\end{table}%
\begin{figure*}[h!]
   \centering
   \begin{subfigure}[t]{0.5\textwidth}\vskip 0pt
       \centering
       \includegraphics[width=3in]{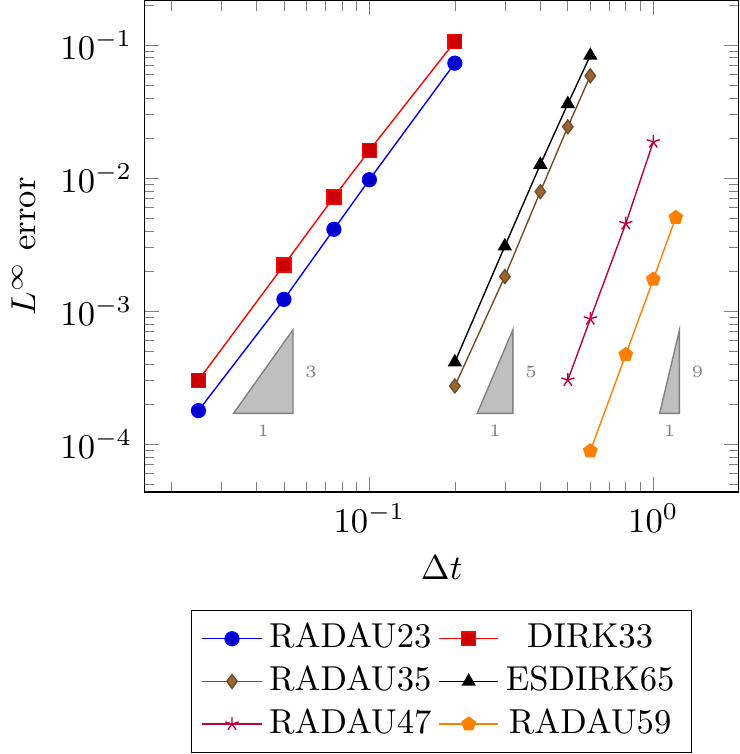}
       \caption{Log-log plot of $L^\infty$ error vs.\ $1/\Delta t$}
       \label{fig:ev-conv}
   \end{subfigure}%
   ~
   \begin{subfigure}[t]{0.5\textwidth}\vskip 0pt
       \centering
       \includegraphics[width=3.12in]{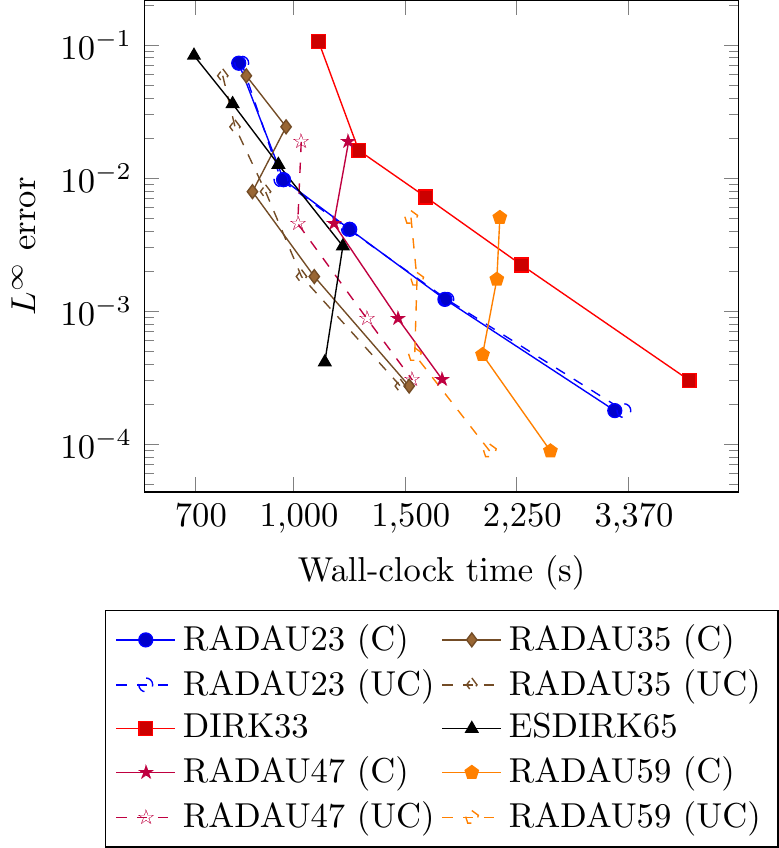}
       \caption{Log-log plot of $L^\infty$ error vs.\ wall-clock time}
       \label{fig:ev-wallclock}
   \end{subfigure}
   \caption{Log-log plots of $L^\infty$ error vs.\ time-step and wall-clock
            time for Euler vortex test case. Stage-coupled preconditioners are
            shown in solid lines, stage-uncoupled in dashed lines.}
\end{figure*}

A log-log plot of the $L^\infty$ error vs.\ $\Delta t$ is shown in Figure
\ref{fig:ev-conv}. A log-log plot of the $L^\infty$ error vs.\ wall-clock time
is shown in Figure \ref{fig:ev-wallclock}. We remark that the RADAU23 method
achieved the same accuracy as the DIRK33 method in faster runtime for all of the
cases considered. Among the fifth-order methods, we did not observe one method
to be clearly more efficient than the others. The difference between the
stage-coupled and stage-uncoupled preconditioners was found to be insignificant
for the third- and fifth-order methods, and the stage-uncoupled preconditioner
resulted in faster performance for the seventh- and ninth-order methods. 

\subsection{High Reynolds number flow over 2D NACA airfoil (LES)}
\label{sec:naca-2d}
\begin{figure}[t!]
   \centering
   \begin{subfigure}[b]{0.5\textwidth}
       \centering
       \includegraphics[width=3in]{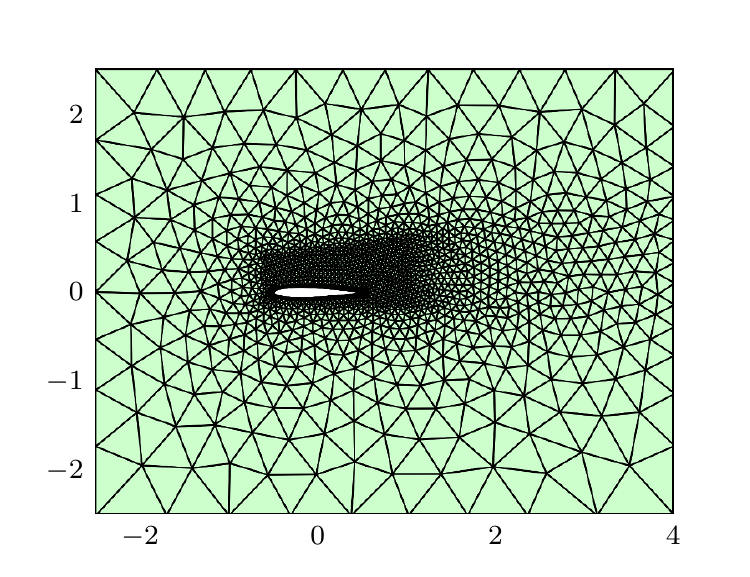}
       \caption{Mesh of NACA airfoil with 3154 elements}
       \label{fig:naca-mesh1}
   \end{subfigure}%
   ~
   \begin{subfigure}[b]{0.5\textwidth}
       \centering
       \includegraphics[width=3in]{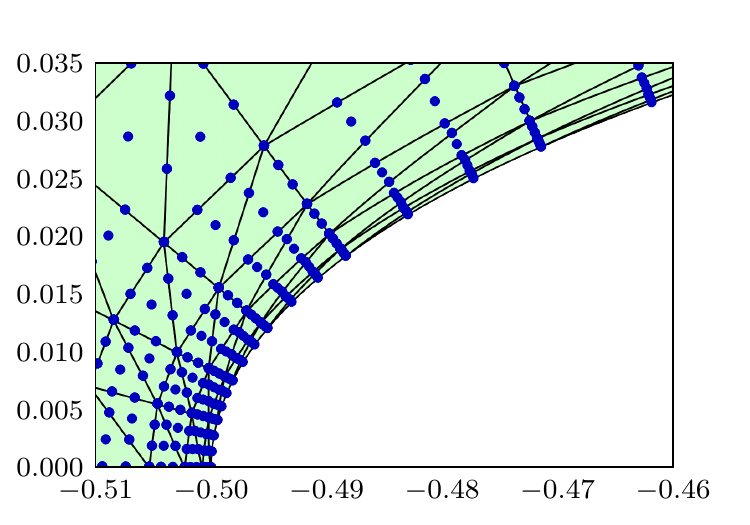}
       \caption{Boundary layer elements with $p=3$ DG nodes}
       \label{fig:naca-bdry}
   \end{subfigure}
   \caption{NACA airfoil mesh}
   \label{fig:naca-mesh}
\end{figure}
In this test case, we consider the two-dimensional viscous flow around a NACA
airfoil with an angle of attack of $\SI{30}{\degree}$. The fluid domain is the
rectangle $[-2.5, 4] \times[-2.5, 2.5]$. For this case, the Mach number is taken
to be 0.1, the gas constant 1.4, and the Reynolds number 40,000.
The wing is centered vertically, and placed closer to the inlet boundary. The
mesh consists of 3154 triangular elements, with finer elements close to the
airfoil and in the wake. The mesh and the boundary-layer elements with DG nodes
are depicted in Figure \ref{fig:naca-mesh}, and are rotated by
$\SI{30}{\degree}$ to correspond with the desired angle of attack. The triangles
near the boundary of the wing are refined in the transverse direction, resulting
in highly anisotropic elements. These stretched elements give rise to a highly
restrictive CFL condition, suggesting that this problem is particularly well
suited to implicit methods. We have found that the CFL condition renders
explicit methods impractical for this problem, with the fourth-order explicit
Runge-Kutta method exhibiting instability for time steps greater than $7\times
10^{-8}\ \si\second$. On the other hand, the implicit Runge-Kutta methods remain
stable for time steps many orders of magnitude larger.

\begin{figure}[t!]
   \centering
   \begin{subfigure}[c]{0.5\textwidth}
       \centering
       \includegraphics[width=3in]{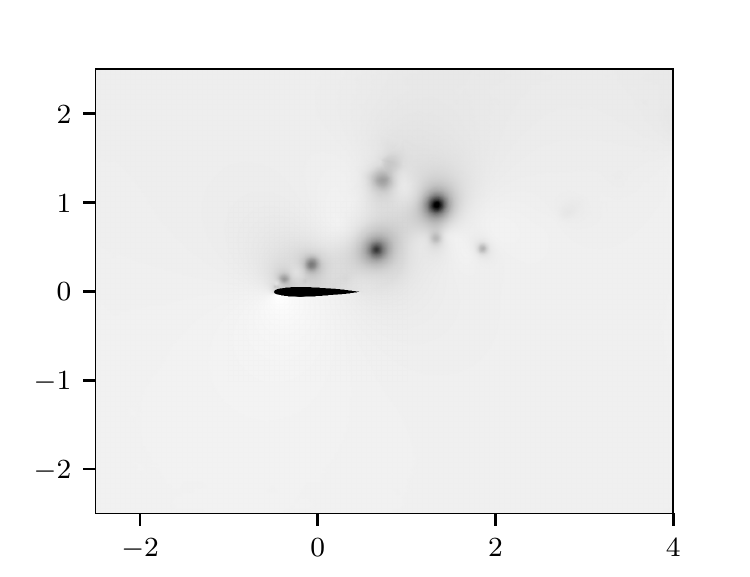}
       \caption{Initial condition, $t=\SI{5}{\s}$ (density)}
       \label{fig:naca-ic}
   \end{subfigure}%
   ~
   \begin{subfigure}[c]{0.5\textwidth}
       \centering
       \includegraphics[width=3in]{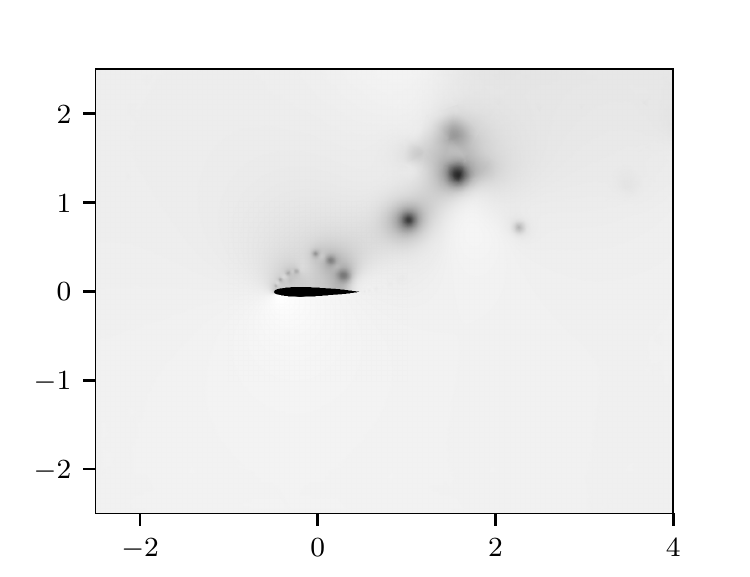}
       \caption{Solution, $t=\SI{5.75}{\s}$ (density)}
       \label{fig:naca-final}

   \end{subfigure}
   \caption{Density plots for NACA $\rm Re = 40k$ test case}
   \label{fig:naca-soln}
\end{figure}

A no-slip condition is enforced on the boundary of the airfoil, and far-field
conditions are enforced on all other boundaries.  At time $t = \SI{0}{\s}$ the
solution is set to be free-stream everywhere. The solution is then integrated
until time $t = \SI{5}{\s}$, at which point vortices have developed in the wake
of the wing. This solution (shown in Figure \ref{fig:naca-ic}) is taken to be
the initial condition for our test case. The finite element space is taken to be
piecewise degree 3 polynomials, with 10 nodes per element, resulting in 
94,620 degrees of freedom.
%

\subsubsection{Solver efficiency}
\begin{figure*}[t!]
   \centering
   \begin{subfigure}[t]{0.5\textwidth}
       \centering
       \includegraphics[width=3in]{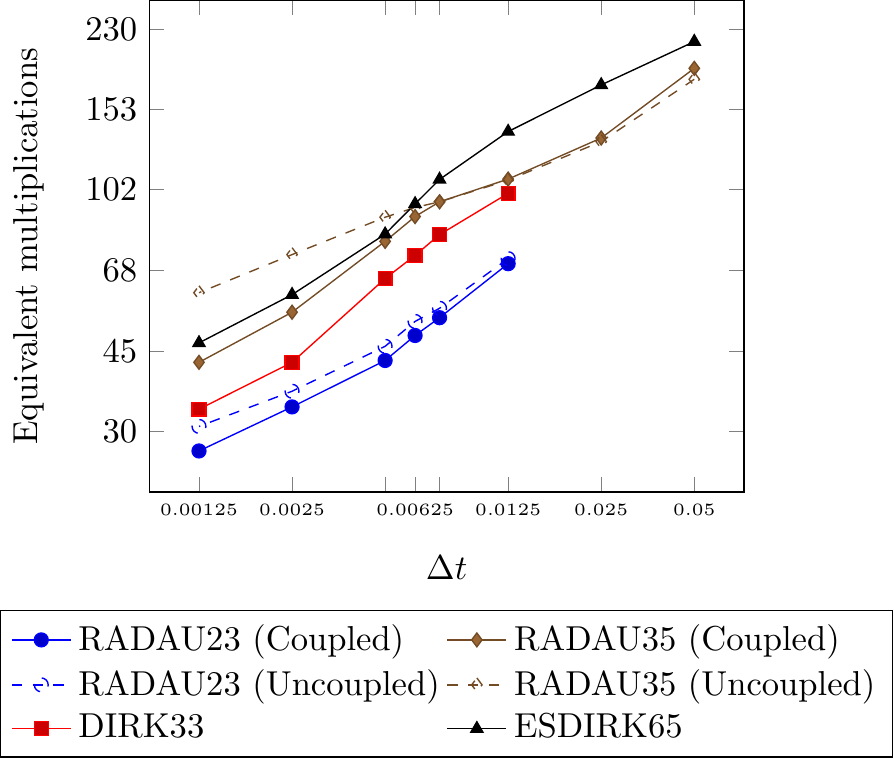}
       \caption{Matrix-vector multiplications vs. $\Delta t$}
       \label{fig:naca-matvec}
   \end{subfigure}%
   ~
   \begin{subfigure}[t]{0.5\textwidth}
       \centering
       \includegraphics[width=3in]{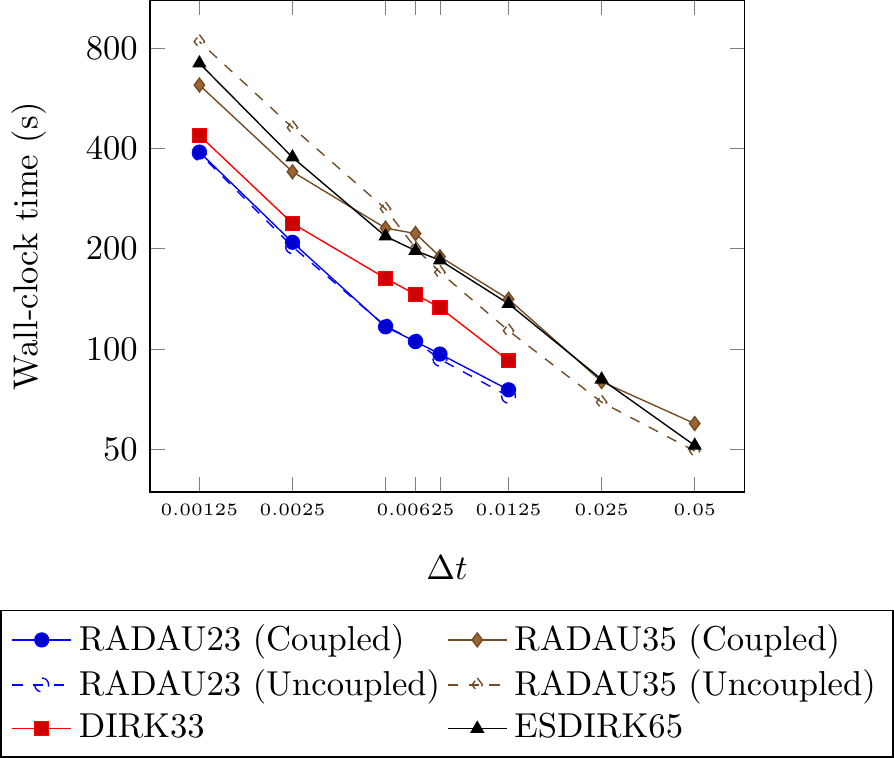}
       \caption{Wall-clock time vs. $\Delta t$}
       \label{fig:naca-wallclock}
   \end{subfigure}
   \caption{Log-log plots of average number of equivalent multiplications 
            and wall-clock time vs. $\Delta t$ for the 2D NACA LES test case}
   \label{fig:naca-results}
\end{figure*}

We study the effectiveness of the DIRK and Radau IIA IRK methods by comparing
both the average number of equivalent multiplications per linear solve and the
total wall-clock time. As in the previous case, the methods RADAU23, RADAU35,
DIRK33, and ESDIRK65 are used. As in the case of the Euler vortex,
a tolerance of $10^{-8}$ is used for the Newton solver. We integrate the
equations from time $t = \SI{5}{\s}$ until $t = \SI{5.75}{\s}$. For the
third-order methods, time steps of $\Delta t =  1.25 \times 10^{-2}, 7.50 \times
10^{-3}, 6.25 \times 10^{-3}, 5.00 \times 10^{-3}, 2.50 \times 10^{-3}, 1.25
\times 10^{-3},$ are used. For the fifth-order methods,  we use the same time
steps, in addition to the larger time steps of $\Delta t = 5.00 \times 10^{-2},$
and $\Delta t = 2.50 \times 10^{-2}$. In Table \ref{tab:naca-results} we present
the runtime and number of equivalent multiplications for all of the methods
considered. As in Section \ref{sec:ev-precond}, we compute the number of
$n\times n$ matrix-vector products performed per Newton iteration (over all the
stages) by multiplying the average number of GMRES iterations by the number of
implicit stages. In Figure \ref{fig:naca-matvec} we present a log-log plot of
the average number of equivalent multiplications vs.\ $\Delta t$. In Figure
\ref{fig:naca-wallclock} we present the total wall-clock time required to run
the simulation until the final time.
\begin{table}[t!]
\scriptsize
\centering
\caption{Equivalent multiplications and wall-clock time for NACA LES test case}
\label{tab:naca-results}
\begin{tabular}{l|cc>{\hspace{1pc}}cc>{\hspace{1pc}}cc}
\toprule
          & \multicolumn{2}{c}{RADAU23 (Coupled)}
          & \multicolumn{2}{c}{RADAU23 (Uncoupled)}
          & \multicolumn{2}{c}{DIRK33} \\
$\Delta t$ & Mult. &Time (s) & Mult. &Time (s) &Mult. &Time (s) \\
\midrule
$1.25\times 10^{-2}$  & 70.2 & 75.4   &	72.0 &	72.1  &   100.2   &   92.1 \\
$7.50\times 10^{-3}$  & 53.4 & 96.5   &	56.0 &	93.1  &   81.3    &   133.0\\
$6.25\times 10^{-3}$  & 48.8 & 105.2  &	52.4 &	105.5 &   73.2    &   145.7\\
$5.00\times 10^{-3}$  & 43.0 & 116.5  &	46.2 &	117.1 &   65.1    &   162.9\\
$2.50\times 10^{-3}$  & 34.0 & 208.9  &	36.8 &	202.7 &   42.6    &   238.0\\
$1.25\times 10^{-3}$  & 27.2 & 389.8  &	30.8 &	386.7 &   33.6    &   437.5\\
\toprule
          & \multicolumn{2}{c}{RADAU35 (Coupled)}
          & \multicolumn{2}{c}{RADAU35 (Uncoupled)}
          & \multicolumn{2}{c}{ESDIRK65} \\
$\Delta t$ &Mult. &Time (s) &Mult. &Time (s) &Mult. &Time (s) \\
\midrule
$5.00\times 10^{-2}$  & 189.2 & 59.7  &	178.5 &	49.4  &   216.0   &   51.3 \\
$2.50\times 10^{-2}$  & 132.6 & 79.8  &	130.2 &	69.2  &   173.5   &   81.0 \\
$1.25\times 10^{-2}$  & 107.7 & 141.0 &	107.1 &	113.4 &   137.0   &   136.6\\
$7.50\times 10^{-3}$  & 96.0  & 189.2 &	96.0  &	169.8 &   107.5   &   184.6\\
$6.25\times 10^{-3}$  & 89.1  & 222.0 &	93.3  &	200.7 &   95.0    &   197.2\\
$5.00\times 10^{-3}$  & 78.6  & 230.4 &	88.8  &	263.2 &   81.5    &   217.8\\
$2.50\times 10^{-3}$  & 54.9  & 340.2 & 73.5  & 461.7 &   60.0    &   375.9\\
$1.25\times 10^{-3}$  & 42.6  & 619.4 & 60.6  & 836.3 &   47.0    &   720.1\\
\bottomrule
\end{tabular}
\end{table}

For the third-order methods, the RADAU23 IRK method required fewer matrix-vector
multiplications than the DIRK33 method, for the same choice of $\Delta t$,
resulting in a shorter run time. The uncoupled preconditioner resulted in a
somewhat larger number of GMRES iterations, and hence more matrix-vector
multiplications, but the difference in run time was found to be negligible. For
the fifth-order methods, the RADAU35 method with stage-coupled preconditioner
resulted in a smaller number of matrix-vector multiplications per linear solve
than the ESDIRK65 method. For larger $\Delta t$, the stage-uncoupled
preconditioner proved to be effective, resulting in faster run times than both
the DIRK method and the stage-coupled preconditioner. For smaller $\Delta t$,
the stage-uncoupled preconditioner required a greater number of GMRES iterations
and hence longer run times.

\subsubsection{Accuracy}
\begin{table}[h!]
\centering
\scriptsize
\caption{$L^\infty$ error for NACA LES test cases}
\label{tab:naca-errors}
\begin{tabular}{l|cc>{\hspace{1pc}}cc>{\hspace{1pc}}c}
\toprule
          & \multicolumn{2}{c}{RADAU23} & \multicolumn{2}{c}{DIRK33} \\
$\Delta t$ &$L^\infty$ & Order & $L^\infty$ & Order & Ratio\\
\midrule
$1.25\times 10^{-2}$  & $6.006\times 10^{-1}$ & -    & $4.514\times 10^{-1}$ & -     & 1.668\\
$7.50\times 10^{-3}$  & $1.743\times 10^{-1}$ & 2.42 & $1.615\times 10^{-1}$ & 2.01  & 1.706\\
$6.25\times 10^{-3}$  & $1.207\times 10^{-1}$ & 2.02 & $1.423\times 10^{-1}$ & 0.70  & 1.717\\
$5.00\times 10^{-3}$  & $8.503\times 10^{-2}$ & 1.57 & $1.216\times 10^{-1}$ & 0.70  & 1.735\\
$2.50\times 10^{-3}$  & $1.347\times 10^{-2}$ & 2.66 & $2.022\times 10^{-2}$ & 2.59  & 1.839\\
$1.25\times 10^{-3}$  & $1.342\times 10^{-3}$ & 3.33 & $2.708\times 10^{-3}$ & 2.90  & 1.823\\
\toprule
          & \multicolumn{2}{c}{RADAU35} & \multicolumn{2}{c}{ESDIRK65} \\
$\Delta t$ &$L^\infty$ & Order & $L^\infty$ & Order & Ratio\\
\midrule
$5.00\times 10^{-2}$  & $4.601\times 10^{-1}$ & -     & $6.239\times 10^{-1}$ & -      & 2.308\\
$2.50\times 10^{-2}$  & $3.158\times 10^{-1}$ & 0.54  & $5.853\times 10^{-1}$ & 0.092  & 1.896\\
$1.25\times 10^{-2}$  & $6.739\times 10^{-2}$ & 2.23  & $1.127\times 10^{-1}$ & 2.377  & 2.414\\
$7.50\times 10^{-3}$  & $5.124\times 10^{-3}$ & 5.04  & $3.989\times 10^{-3}$ & 6.540  & 4.441\\
$6.25\times 10^{-3}$  & $5.329\times 10^{-3}$ & -0.22 & $2.147\times 10^{-3}$ & 3.397  & 3.263\\
$5.00\times 10^{-3}$  & $2.616\times 10^{-3}$ & 3.19  & $3.279\times 10^{-3}$ & -1.897 & 3.200\\
$2.50\times 10^{-3}$  & $8.135\times 10^{-4}$ & 1.68  & $2.222\times 10^{-3}$ & 0.561  & 1.044\\
$1.25\times 10^{-3}$  & $2.059\times 10^{-4}$ & 1.98  & $1.727\times 10^{-4}$ & 3.686  & 0.896\\
\bottomrule
\end{tabular}
\end{table}
\begin{figure*}[h!]
   \centering
   \begin{subfigure}[t]{0.5\textwidth}
       \centering
       \includegraphics[width=3in]{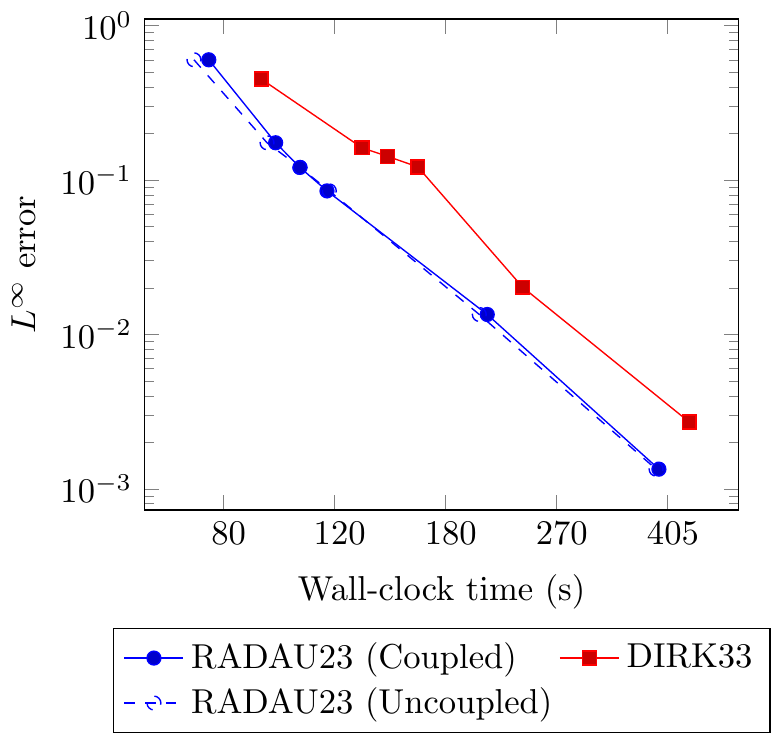}
   \end{subfigure}%
   ~
   \begin{subfigure}[t]{0.5\textwidth}
       \centering
       \includegraphics[width=3in]{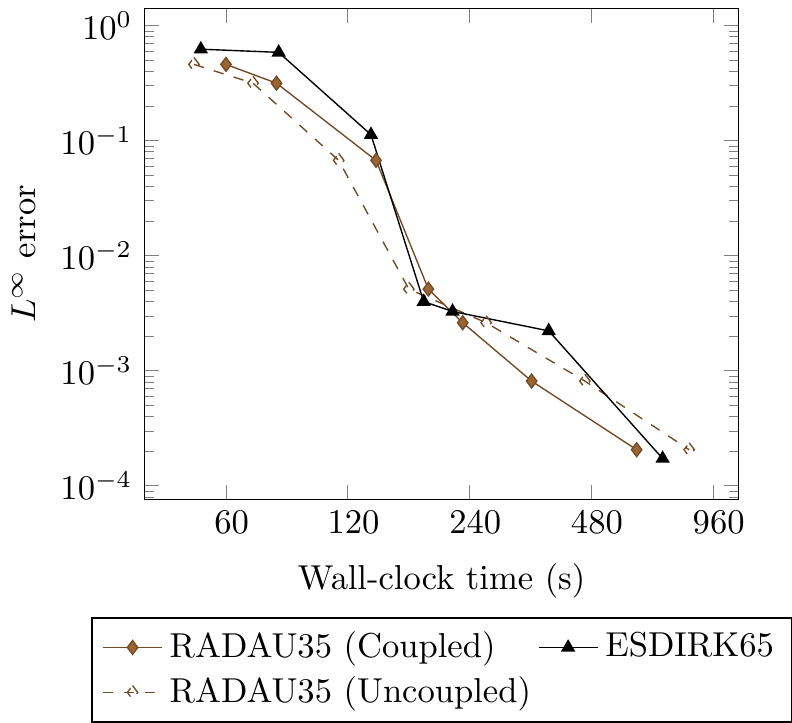}
   \end{subfigure}
   \caption{Log-log plots of $L^\infty$ errors vs.\ wall-clock time for the
            NACA LES test case}
   \label{fig:naca-errors}
\end{figure*}

We note that the above comparisons were made for equal choices of $\Delta t$.
Because the Radau IIA method enjoy a smaller leading coefficient of the
truncation error, we expect to achieve better accuracy for the same time step.
We therefore study the accuracy of the methods applied to the above problem by
considering the semidiscrete system of equations purely as a system of ODEs. We
compute a reference solution by numerically integrating the equations for 6000
time steps with $\Delta t = 1.25 \times 10^{-4}$ using the fifth-order ESDIRK
method. Then, we take this solution to be the ``exact'' solution, with respect
to which the $L^\infty$ norm of the error is computed, and perform
a grid convergence study. For each choice of time step, we compute the
$L^\infty$ norm of the error, the approximate rate of convergence of the method,
and the ratio of the DIRK error to the Radau IIA error in Table
\ref{tab:naca-errors}. In Figure \ref{fig:naca-errors}, we present log-log plots
of the wall-clock time vs.\ $L^\infty$ error. We note that we do not observe the
formal order of temporal accuracy for this test problem, possibly due to the
choice of time step, which is about five orders of magnitude larger than the 
explicit CFL, together with the stiff and turbulent nature of the problem.
In order to verify the formal order of accuracy, we also run this test problem
with time steps that are two to three orders of magnitude smaller than those 
of the above comparison (for a shorted total time of integration). The results 
presented in Table \ref{tab:naca-small-dt} confirm that the expected theoretical
orders of accuracy are attained for all methods considered.

\begin{table}[h!]
\centering
\scriptsize
\caption{$L^\infty$ error for NACA LES test cases (order verification)}
\label{tab:naca-small-dt}
\begin{tabular}{l|cc>{\hspace{1pc}}cc>{\hspace{1pc}}c}
\toprule
          & \multicolumn{2}{c}{RADAU23} & \multicolumn{2}{c}{DIRK33} \\
$\Delta t$ &$L^\infty$ & Order & $L^\infty$ & Order \\
\midrule
$6.25\times 10^{-5}$    & $5.694\times 10^{-6}$ & -    & $3.580\times 10^{-6}$ & -    \\
$3.125\times 10^{-5}$   & $8.765\times 10^{-7}$ & 2.70 & $5.004\times 10^{-7}$ & 2.84 \\
$1.5625\times 10^{-5}$  & $1.177\times 10^{-7}$ & 2.90 & $6.401\times 10^{-8}$ & 2.97 \\
$7.8125\times 10^{-6}$  & $1.494\times 10^{-8}$ & 2.98 & $8.033\times 10^{-9}$ & 2.99 \\
\toprule
          & \multicolumn{2}{c}{RADAU35} & \multicolumn{2}{c}{ESDIRK65} \\
$\Delta t$ &$L^\infty$ & Order & $L^\infty$ & Order \\
\midrule
$1.25\times 10^{-4}$    & $3.438\times 10^{-6}$ & -     & $1.231\times 10^{-6}$ & -       \\
$6.25\times 10^{-5}$    & $1.621\times 10^{-7}$ & 4.41  & $4.758\times 10^{-8}$ & 4.69    \\
$3.125\times 10^{-5}$   & $5.434\times 10^{-9}$ & 4.90  & $1.602\times 10^{-9}$ & 4.89    \\
$1.5625\times 10^{-5}$  & $1.964\times 10^{-10}$ & 4.79  & $7.913\times 10^{-11}$ & 4.34  \\
\bottomrule
\end{tabular}
\end{table}

We remark that for the third-order methods, we can achieve the same accuracy as
the DIRK33 method with a faster run time using the RADAU23 method, with both
the stage-coupled or stage-uncoupled preconditioner. The differences in run time
between the two preconditioners were negligible. For the fifth-order methods,
our results indicate that the RADAU35 method outperformed the ESDIRK for a 
majority of the test cases considered. For this test problem, we found the
stage-uncoupled preconditioner to perform better for the larger choices of
$\Delta t$.

\subsubsection{Parallel performance on NACA airfoil}
\label{sec:parallel-results}
We study the parallel performance of the IRK and DIRK solvers applied to the
two-dimensional NACA airfoil. We choose a representative time step of $\Delta t
= 1.25\times 10^{-2}$, and integrate the system for five time steps until $t =
\SI{5.0625}{\s}$, using as before the solution at $t = \SI{5}{\s}$ as the
initial condition. We perform this test using both the Radau IIA IRK and the
DIRK solvers. For the Radau methods, we use both the stage-coupled and
stage-uncoupled ILU(0) preconditioners.

\begin{figure}[b!]
   \centering
   \includegraphics[width=5.8in]{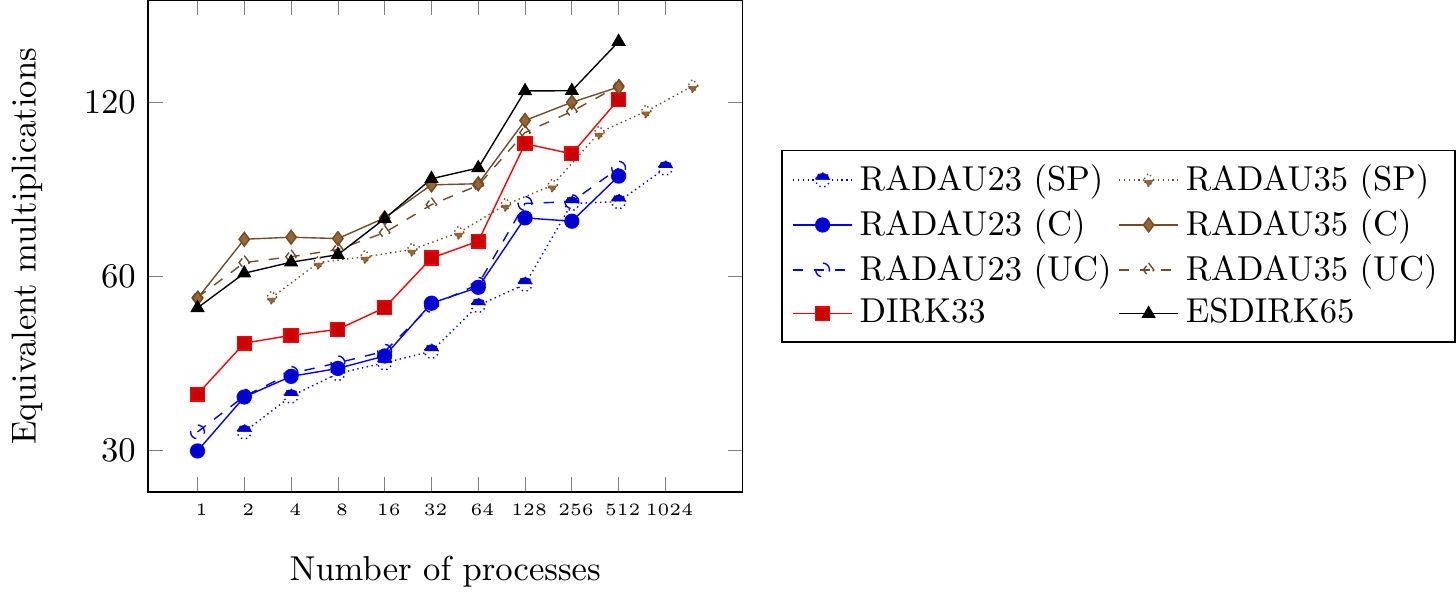}
   \caption{Log-log plot of average number of equivalent multiplications
            vs.\ number of processes. Dashed lines indicate stage-uncoupled
            block ILU(0) preconditioner, and dotted lines indicate
            stage-parallel solver.}
   \label{fig:naca-parallel}
\end{figure}

Using the method described in Section \ref{sec:parallel}, we decompose the
domain into a set number of partitions according to the number of processes.
For the DIRK and stage-coupled IRK solvers, the number of partitions is equal
to the number of processes. For the stage-uncoupled IRK solver, we can choose
the number of partitions to be a factor of $s$ smaller than the number of
processes. We consider the mesh decomposed into 1, 2, 4, 8, 16, 32, 64, 128,
256, and 512 partitions. Then, we compute the average number of GMRES iterations
required per solve. In the case of the DIRK methods, we multiply the number
of iterations by the number of implicit stages to obtain the number of
equivalent multiplications performed assuming one Newton iteration.
Similarly, in the case of the Radau IIA IRK methods, we multiply the number of
iterations by the number of stages to obtain the number of $n \times n$
matrix-vector multiplications. In Figure \ref{fig:naca-parallel} we show
a log-log plot of the average number of equivalent multiplications vs.\
number of processes.

As the number of processes (and hence number of mesh partitions) increases, we
observe an increase in the number of GMRES iterations required to converge. This
is because the contributions between different mesh partitions are ignored in
the block ILU(0) factorization, rendering the preconditioner less effective.
Since the stage-uncoupled ILU(0) preconditioner can be parallelized for the same
number of processes with a factor of $s$ fewer partitions, we
notice that approximately 15--20\% fewer matrix-vector multiplications are
required when compared with the standard uncoupled solver. This suggests a
substantial benefit to the stage-uncoupled block ILU(0) preconditioner when run
in a massively-parallel environment. This difference in performance is
numerically validated in the following three-dimensional NACA test case.

\subsection{Parallel large eddy simulation (LES) of 3D NACA airfoil}
\label{sec:naca-3d}
\begin{figure}[h!]
  \centering
  \begin{subfigure}[t]{0.5\textwidth}
      \centering
{%
\setlength{\fboxsep}{0pt}%
\setlength{\fboxrule}{0.5pt}%
\fbox{\includegraphics[width=3in]{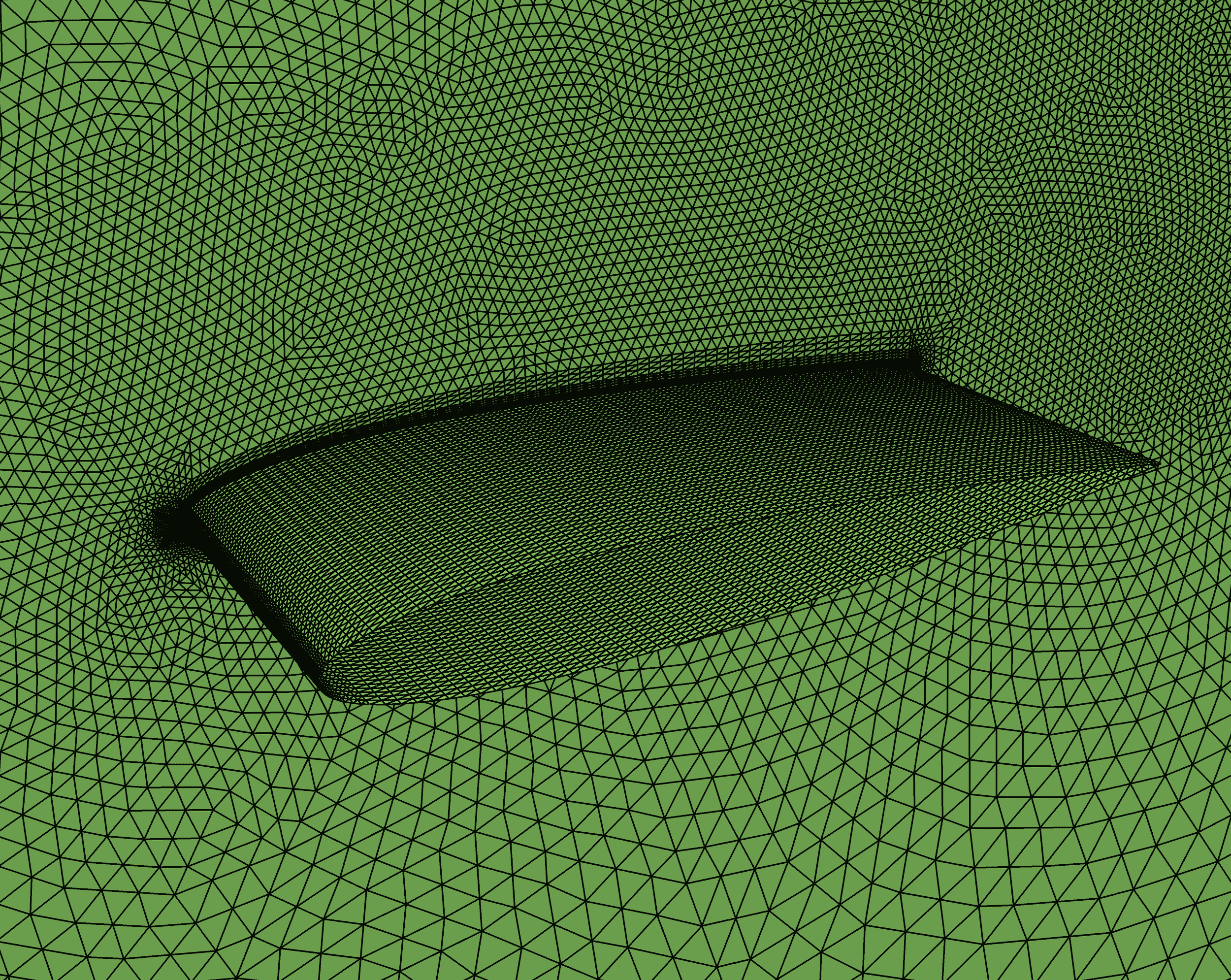}}
}%
      \caption{Boundaries of three-dimensional NACA mesh.}
      \label{fig:naca3d-mesh}
  \end{subfigure}%
  ~
  \begin{subfigure}[t]{0.5\textwidth}
      \centering
{%
\setlength{\fboxsep}{0pt}%
\setlength{\fboxrule}{0.5pt}%
\includegraphics[width=3in]{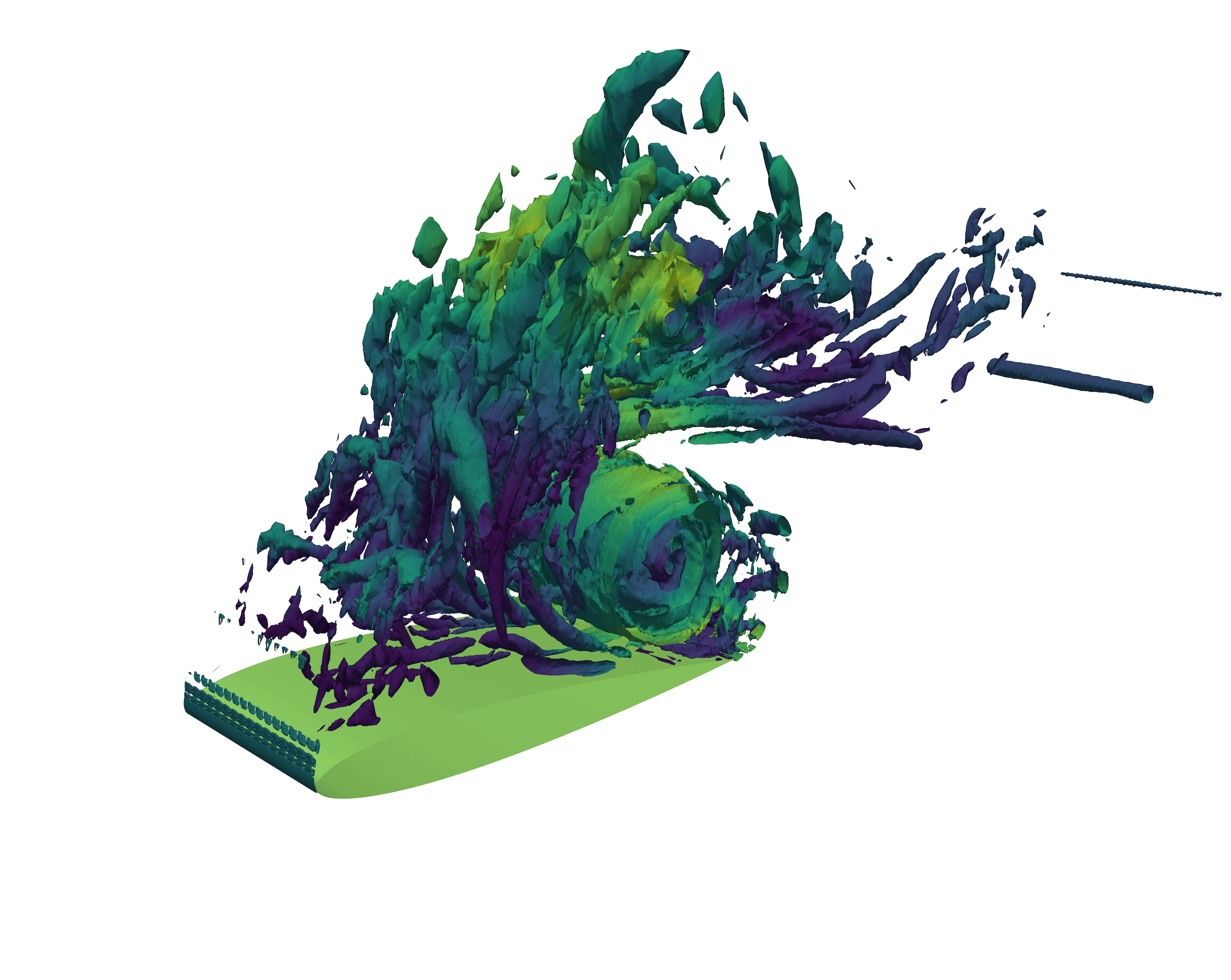}
}%
      \caption{Isosurfaces of $Q$-criterion, $Q = 25$, colored by velocity 
               magnitude.}
      \label{fig:naca3d-isosurface}
  \end{subfigure}
  \caption{Three-dimensional NACA LES test case.}
  \label{fig:naca3d}
\end{figure}
For our final test case, we consider the three-dimensional viscous flow over a
NACA airfoil with angle of attack of $\SI{30}{\degree}$. The Reynolds number is
taken to be 5,300 and the Mach number 0.1. The governing equations
are given by equations \eqref{eq:ns-1} through \eqref{eq:ns-3} with the
isentropic assumption discussed in Section \ref{sec:eqns}. The mesh consists of
151,392 tetrahedral elements. The local basis consists of degree 3 polynomials,
for a total of 20 nodes per element. We consider the fluid to be isentropic, and
hence the solution consists of 4 components, resulting in a total of 
12,111,360
degrees of freedom. The mesh is shown in Figure \ref{fig:naca3d-mesh}.

A no-slip condition is enforced on the boundary
of the airfoil, periodic conditions are enforced in the span-wise direction,
and far-field conditions are enforced on all other boundaries. 
The solution is initialized to freestream conditions, and then integrated
numerically until time $t = \SI{4}{\s}$ using $\Delta t = \SI{0.02}{\s}$. At
this point, vortices have developed in the wake of the wing. In Figure
\ref{fig:naca3d-isosurface}, isosurfaces of the $Q$-criterion, for $Q = 25$, are
shown. The $Q$-criterion, proposed by Dubief and Delcayre in \cite{Dubief2000}
is often used to identify vortical structures, and is defined as the difference
of the symmetric and antisymmetric components of the velocity gradient,
\begin{equation}
  Q = \frac{1}{2}\left( \Omega_{ij}\Omega_{ij} - S_{ij}S_{ij} \right),
\end{equation}
where $\Omega_{ij} = \frac{1}{2}\left(u_{i,j} - u_{j,i}\right)$, and 
$S_{ij} = \frac{1}{2}\left(u_{i,j} + u_{j,i}\right)$.

\begin{figure}[h]
   \centering
   \includegraphics[width=5.8in]{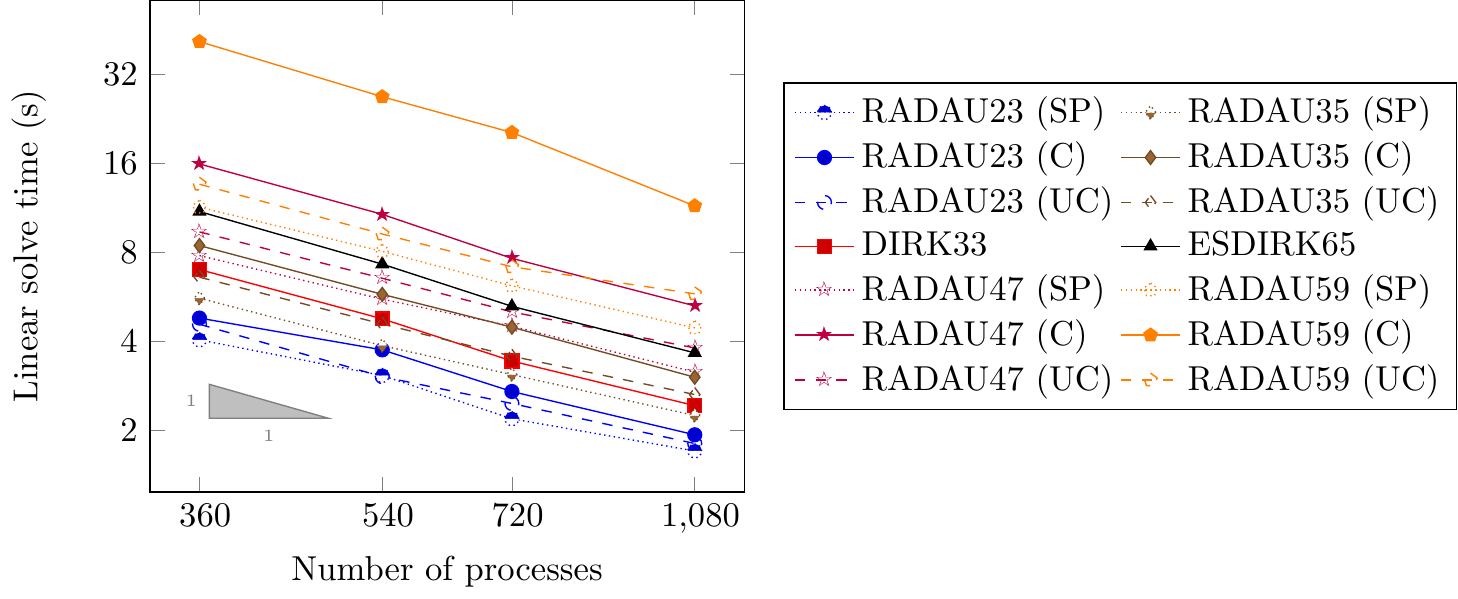}
   \caption{Log-log plot of linear solve time in seconds
            vs.\ number of processes. Dashed lines indicate stage-uncoupled
            block ILU(0) preconditioner, and dotted lines indicate
            stage-parallel solver. A reference triangle for perfect speedup
            is shown.}
   \label{fig:naca3d-runtime}
\end{figure}

We then integrate the equations for 30 time
steps using the third- and fifth-order DIRK methods, as well as the Radau IIA
methods of order 3, 5, 7, and 9. Due to the turbulent nature of this problem, we
do not study the accuracy of the numerical solutions, but rather the efficiency
of the solvers for fixed $\Delta t$. We also consider the parallel scaling of
the solvers by running our test case on 360, 540, 720, and 1080 processes. For
the Radau IIA methods, we consider three preconditioners: stage-coupled,
stage-uncoupled, and stage-parallel. As with the DIRK methods, for the
stage-coupled and stage-uncoupled preconditioners we decompose the mesh into a
number of partitions equal to the number of processes. For the stage-parallel
ILU(0) preconditioner we decompose the mesh into a factor of $s$ fewer
partitions. The preconditioner then exploits the stage parallelism using the
methodology described in Section \ref{sec:parallel}. We record the average
wall-clock time required per linear solve using each of the methods, and present
the results in Figure \ref{fig:naca3d-runtime}.

For each order of accuracy, the Radau IIA method with the stage-parallel
preconditioner resulted in the fastest runtime. For orders three and five, for
which we compare against equal-order DIRK methods, the Radau IIA methods
resulted in faster performance with all of the preconditioners considered. Since
the number of stages is greater for the higher order methods, and the
stage-parallel preconditioner allows for a factor of $s$ fewer mesh partitions,
we would expect that the reduction in the number of GMRES iterations due to the
improved ILU preconditioner would be greater than for the lower order methods.
Indeed, we observe that the relative performance gain of the stage-parallel
preconditioner compared with the stage-uncoupled preconditioner increases as the
the order of the method increases. For example, while the
third-order, stage-parallel preconditioner is only about 5--10\% faster than the
stage-uncoupled preconditioner, and about 20\% faster than the stage-coupled
preconditioner, the ninth-order stage-parallel preconditioner is between
20--30\% faster than the stage-uncoupled preconditioner, and approximately twice
as fast as the stage-coupled preconditioner. Furthermore, as the number of
processes approaches the strong scaling limit, we anticipate that the additional
factor of $s$ processes allowed by the stage-parallel preconditioner will result
in even greater benefit, especially for those methods with a large number of
stages.

\section{Conclusion}
In this paper, we have developed a new strategy for efficiently solving the
large, coupled linear systems arising from fully implicit Runge-Kutta time
discretizations. By transforming the system of equations, the computational work
required per GMRES iterations is reduced significantly. This new method makes it
feasible to use high-order, $L$-stable implicit Runge-Kutta methods, such as the
Radau IIA methods, as time integrators for discontinuous Galerkin
discretizations. We additionally develop new preconditioners for these methods,
included a parallel-in-time ILU preconditioner that allows for the computation 
of the Runge-Kutta stages simultaneously.

Numerical experiments on both two- and three-dimensional fluid flow problems are
performed using the Radau IIA methods of up to ninth order. These results
indicate that using the transformed system of equations, the fully implicit IRK
methods are competitive with, and in our experience, often preferable to the
more standard DIRK methods, both in terms of efficiency and accuracy. In a
parallel computing environment, the stage-parallel ILU preconditioner results in
additional performance gains.

\section{Acknowledgments}
This research used resources of the National Energy Research Scientific 
Computing Center, a DOE Office of Science User Facility supported by the Office 
of Science of the U.S. Department of Energy under Contract No.\ 
DE-AC02-05CH11231. The first author was supported by the Department of Defense 
through the National Defense Science \& Engineering Graduate Fellowship Program.

\appendix
\input{tableaux}

\bibliographystyle{plain}
\bibliography{IRK}

\end{document}

%% file: tableaux.tex

\newgeometry{margin=2cm} 
\begin{landscape}
\setlength\extrarowheight{6pt}
\newcommand*{\MmaComment}[1]{\hfill\makebox[10.0cm][l]{\rmfamily #1}}%

\section{Butcher tableaux} \label{app:tableaux}

\begin{tabular}{llll}
RADAU23: &
$\begin{aligned}
   \begin{array}{r|rr}
       1/3  & 5/12 & -1/12 \\
       1    & 3/4  & 1/4   \\
       \hline
            & 3/4  & 1/4
   \end{array}
\end{aligned}$ &
\hspace{1in}RADAU35: &
$\begin{aligned}
   \begin{array}{r|rrr}
       \tfrac{2}{5} - \tfrac{\sqrt{6}}{10} 
           & \tfrac{11}{45} - \tfrac{7\sqrt{6}}{360}
           & \tfrac{37}{225} - \tfrac{169\sqrt{6}}{1800}
           & -\tfrac{2}{225} + \tfrac{\sqrt{6}}{75} \\
       \tfrac{2}{5} + \tfrac{\sqrt{6}}{10} 
           & \tfrac{37}{225} - \tfrac{169\sqrt{6}}{1800}
           & \tfrac{11}{45} - \tfrac{7\sqrt{6}}{360}
           & -\tfrac{2}{225} - \tfrac{\sqrt{6}}{75} \\
       1
           & \tfrac{4}{9} - \tfrac{\sqrt{6}}{36}
           & \tfrac{4}{9} + \tfrac{\sqrt{6}}{36}
           & \tfrac{1}{9}  \vphantom{\dfrac{1}{2}} \\
       \hline
           & \tfrac{4}{9} - \tfrac{\sqrt{6}}{36}
           & \tfrac{4}{9} + \tfrac{\sqrt{6}}{36}
           & \tfrac{1}{9}
   \end{array}
\end{aligned}$ \\[1in]
\multicolumn{4}{l}{}
\varwidth{\linewidth}
The remaining Radau IIA (4 and 5 stage) tableaux were generated using the 
following Mathematica code, according to the derivation in \cite{Axelsson1969}, 
where \texttt{c[n]}, and \texttt{A[n]} give the abscissa, and Butcher matrix, 
respectively, for the $n$-stage Radau IIA method.
\begin{lstlisting}
   q[n_, x_] = LegendreP[n, 2 x - 1] - LegendreP[n - 1, 2 x - 1]; !
    \MmaComment{Define the polynomial $Q_n(x)$ according to equation (2.1) %
    from \cite{Axelsson1969}}!
   c[n_] := (x /. NSolve[q[n, x] == 0, x]) // Re // (Sort[#, Less] &) !
    \MmaComment{The abscissa $c_k$ are given by the zeros of $Q_n(x)$} !
   !$\ell$![n_, k_, x_] := q[n, x]/((x - c[n][[k]]) (D[q[n, x], x] /. x -> c[n][[k]]))!
    \MmaComment{Define $\ell_k(x)$ as in \cite{Axelsson1969}} !
   a[n_, i_, k_] := NIntegrate[!$\ell$![n, k, x], {x, 0, c[n][[i]]}] !
    \MmaComment{The quadrature coefficients $a_{ik}$ are computed as %
    $\int_0^{c_i} \ell_k(x)~dx$} !
   A[n_] := Table[a[n, i, k], {i, 1, n}, {k, 1, n}] !
    \MmaComment{Finally, we form the Butcher matrix $A$} !
\end{lstlisting}
\endvarwidth\hfill
\\[1in]
DIRK33: &
\multicolumn{3}{l}{
$\begin{aligned}
  \begin{array}{r|rrr}
    \alpha & \alpha          & 0      & 0\\
    \tau_2 & \tau_2 - \alpha & \alpha & 0 \\
    1      & b_1             & b_2    & \alpha\\
    \hline
        & b_1 & b_2 & \alpha
  \end{array}
  \hspace{1in}
  \begin{array}{l}
     \alpha = 1+\frac{\sqrt{6}}{2}  \sin\left(\frac{1}{3}
            \arctan\left(\frac{\sqrt{2}}{4}\right)\right)
         - \frac{\sqrt{2}}{2} \cos\left(\frac{1}{3}
            \arctan\left(\frac{\sqrt{2}}{4}\right)\right)\\
     \tau_2 = (1+\alpha)/2\\
     b_1 = -(6\alpha^2 - 16\alpha +1)/4\\
     b_2 = (6\alpha^2 - 20\alpha +5)/4
  \end{array}
\end{aligned}$
}\\[1in]
ESDIRK65: &
\multicolumn{3}{l}{
$\begin{aligned} \footnotesize
  \begin{array}{r|rrrrrr}
      0
          &  0 & 0 & 0 & 0 & 0 & 0 \\
      0.556107682272893
          &   0.2780538411364465 &  0.2780538411364465 &  0 &  0 &  0 &  0 \\
      1.028127096688746
          &   0.3137405401502951 &  0.4363327154020044 &  0.2780538411364465 &  0 &  0 &  0 \\
      0.540645375074761
          &   0.2741986534107860 &  -0.0164268277321164 &  0.0048197082596452 & 
              0.2780538411364465 &  0 &  0 \\
      0.058741042826253
          &   -0.2441776975175844 &  -3.3203529439447852 &  0.0477747285706825 & 
              3.2974431145814931 &  0.2780538411364465 &  0 \\
      1
          &   -0.2786732780227907 &  1.8929947094010862 &  -0.1280948204262490 & 
              -1.3574693381380240 &  0.5931888860495311 &  0.2780538411364465 \\
      \hline
          & -0.2786732780227907 &  1.8929947094010862 &  -0.1280948204262490 & 
            -1.3574693381380240 &  0.5931888860495311 &  0.2780538411364465
  \end{array}
\end{aligned}$
} 
\end{tabular}

\end{landscape}
\restoregeometry

%% file: IRK.bbl
\begin{thebibliography}{10}

\bibitem{Alexander:1977dk}
Roger Alexander.
\newblock {Diagonally implicit Runge-Kutta methods for stiff O.D.E.'s}.
\newblock {\em SIAM Journal on Numerical Analysis}, 14(6):1006--1021, 1977.

\bibitem{Axelsson1969}
Owe Axelsson.
\newblock A class of {$A$}-stable methods.
\newblock {\em BIT Numerical Mathematics}, 9(3):185--199, 1969.

\bibitem{Bassi2015}
F.~Bassi, L.~Botti, A.~Colombo, A.~Ghidoni, and F.~Massa.
\newblock Linearly implicit {R}osenbrock-type {R}unge-{K}utta schemes applied
  to the discontinuous {G}alerkin solution of compressible and incompressible
  unsteady flows.
\newblock {\em Computers \& Fluids}, 118:305--320, 2015.

\bibitem{Bijl:2002}
Hester Bijl, Mark~H. Carpenter, Veer~N. Vatsa, and Christopher~A. Kennedy.
\newblock Implicit time integration schemes for the unsteady compressible
  {N}avier-{S}tokes equations: Laminar flow.
\newblock {\em Journal of Computational Physics}, 179(1):313--329, 2002.

\bibitem{Boom:2013}
Pieter~D. Boom and David~W. Zingg.
\newblock High-order implicit temporal integration for unsteady compressible
  fluid flow simulation.
\newblock In {\em 21st {AIAA} Computational Fluid Dynamics Conference}.
  American Institute of Aeronautics and Astronautics ({AIAA}), 2013.

\bibitem{Burrage1999}
Kevin Burrage, Craig Eldershaw, Roger Sidje, et~al.
\newblock A parallel matrix-free implementation of a {R}unge-{K}utta code.
\newblock In {\em Joint Australian-Taiwanese Workshop on Analysis and
  Applications}, pages 83--88. Centre for Mathematics and its Applications,
  Mathematical Sciences Institute, The Australian National University, 1999.

\bibitem{Butcher:1976tt}
J.~C. Butcher.
\newblock {On the implementation of implicit Runge-Kutta methods}.
\newblock {\em BIT Numerical Mathematics}, 16(3):237--240, 1976.

\bibitem{Carpenter:2005}
M.~H. Carpenter, C.~A. Kennedy, Hester Bijl, S.~A. Viken, and Veer~N. Vatsa.
\newblock Fourth-order {R}unge-{K}utta schemes for fluid mechanics
  applications.
\newblock {\em Journal of Scientific Computing}, 25(1):157--194, 2005.

\bibitem{Carpenter2003:ef}
Mark~H. Carpenter, Sally~A. Viken, and Eric~J. Nielsen.
\newblock The efficiency of high order temporal schemes.
\newblock {\em AIAA Paper}, 86:2003, 2003.

\bibitem{Dahlquist:1963}
Germund~G. Dahlquist.
\newblock A special stability problem for linear multistep methods.
\newblock {\em BIT Numerical Mathematics}, 3(1):27--43, 1963.

\bibitem{DeSwart1998}
Jacques~JB De~Swart, Walter~M Lioen, and Wolter~A Van Der~Veen.
\newblock {\em Specification of PSIDE}.
\newblock Stichting Mathematisch Centrum, 1998.

\bibitem{Dubief2000}
Yves Dubief and Franck Delcayre.
\newblock On coherent-vortex identification in turbulence.
\newblock {\em Journal of Turbulence}, 1(1):011--011, 2000.

\bibitem{Frank:1985kv}
Reinhard Frank, Josef Schneid, and Christoph~W Ueberhuber.
\newblock {Order Results for Implicit Runge{\textendash}Kutta Methods Applied
  to Stiff Systems}.
\newblock {\em SIAM Journal on Numerical Analysis}, 22(3):515--534, 1985.

\bibitem{Hairer:2013vl}
Ernst Hairer and Gerhard Wanner.
\newblock {B-Convergence}.
\newblock In {\em Solving Ordinary Differential Equations II}, pages 225--238.
  Springer Berlin Heidelberg, Berlin, Heidelberg, 1996.

\bibitem{Hairer:2013tj}
Ernst Hairer and Gerhard Wanner.
\newblock {Construction of Implicit Runge-Kutta Methods}.
\newblock In {\em Solving Ordinary Differential Equations II}, pages 71--90.
  Springer Berlin Heidelberg, Berlin, Heidelberg, 1996.

\bibitem{Jay1999}
Laurent~O. Jay and Thierry Braconnier.
\newblock A parallelizable preconditioner for the iterative solution of
  implicit {R}unge–{K}utta-type methods.
\newblock {\em Journal of Computational and Applied Mathematics}, 111(1–2):63
  -- 76, 1999.

\bibitem{Kanner2015}
Samuel Kanner and Per-Olof Persson.
\newblock Validation of a high-order large-eddy simulation solver using a
  vertical-axis wind turbine.
\newblock {\em AIAA Journal}, 54(1):101--112, 2015.

\bibitem{Kvaerno:2004}
A.~Kv{\ae}rn{\o}.
\newblock Singly diagonally implicit {R}unge-{K}utta methods with an explicit
  first stage.
\newblock {\em BIT Numerical Mathematics}, 44(3):489--502, 2004.

\bibitem{Lambert:1991vi}
J.~D. Lambert.
\newblock {\em {Numerical methods for ordinary differential systems}}.
\newblock John Wiley {\&} Sons, Ltd., Chichester, 1991.

\bibitem{Manteuffel1980}
Thomas~A Manteuffel.
\newblock An incomplete factorization technique for positive definite linear
  systems.
\newblock {\em Mathematics of Computation}, 34(150):473--497, 1980.

\bibitem{Nigro2014MEBDF}
A.~Nigro, A.~Ghidoni, S.~Rebay, and F.~Bassi.
\newblock Modified extended {BDF} scheme for the discontinuous {G}alerkin
  solution of unsteady compressible flows.
\newblock {\em International Journal for Numerical Methods in Fluids},
  76(9):549--574, 2014.

\bibitem{Nigro2014TIAS}
Alessandra Nigro, Carmine~De Bartolo, Francesco Bassi, and Antonio Ghidoni.
\newblock Up to sixth-order accurate {$A$}-stable implicit schemes applied to
  the discontinuous {G}alerkin discretized {N}avier–{S}tokes equations.
\newblock {\em Journal of Computational Physics}, 276:136 -- 162, 2014.

\bibitem{Peraire:2008}
Jaime Peraire and Per-Olof Persson.
\newblock The compact discontinuous {G}alerkin ({CDG}) method for elliptic
  problems.
\newblock {\em SIAM Journal on Scientific Computing}, 30(4):1806--1824, 2008.

\bibitem{Peraire:2011}
Jaime Peraire and Per-Olof Persson.
\newblock High-order discontinuous {G}alerkin methods for {CFD}.
\newblock In Z.~J. Wang, editor, {\em Adaptive High-Order Methods in Fluid
  Dynamics}, chapter~5, pages 119--152. World Scientific, 2011.

\bibitem{Persson2009scalable}
Per-Olof Persson.
\newblock Scalable parallel {N}ewton-{K}rylov solvers for discontinuous
  {G}alerkin discretizations.
\newblock In {\em Proceedings of the 47th AIAA Aerospace Sciences Meeting and
  Exhibit}, 2009.

\bibitem{Persson:2011LES}
Per-Olof Persson.
\newblock High-order {LES} simulations using implicit-explicit {R}unge-{K}utta
  schemes.
\newblock In {\em Proceedings of the 49th AIAA Aerospace Sciences Meeting and
  Exhibit, AIAA}, volume 684, 2011.

\bibitem{Persson:2006}
Per-Olof Persson and Jaime Peraire.
\newblock An efficient low memory implicit {DG} algorithm for time dependent
  problems.
\newblock {\em AIAA paper}, 113:2006, 2006.

\bibitem{Persson:2008by}
Per-Olof Persson and Jaime Peraire.
\newblock {N}ewton-{GMRES} preconditioning for discontinuous {G}alerkin
  discretizations of the {N}avier-{S}tokes equations.
\newblock {\em SIAM Journal on Scientific Computing}, 30(6):2709--2733, 2008.

\bibitem{Reed:1973}
W.~H. Reed and T.~R. Hill.
\newblock Triangular mesh methods for the neutron transport equation.
\newblock {\em Los Alamos Report LA-UR-73-479}, 1973.

\bibitem{Wang:2013gj}
Z.~J. Wang, Krzysztof Fidkowski, R{\'e}mi Abgrall, Francesco Bassi, Doru
  Caraeni, Andrew Cary, Herman Deconinck, Ralf Hartmann, Koen Hillewaert, H.~T.
  Huynh, Norbert Kroll, Georg May, Per-Olof Persson, Bram Leer, and Miguel
  Visbal.
\newblock {High-order CFD methods: current status and perspective}.
\newblock {\em International Journal for Numerical Methods in Fluids},
  72(8):811--845, July 2013.

\end{thebibliography}
